\documentclass[a4paper,10pt]{article}
\usepackage[top=2.5cm,bottom=2.5cm,left=2.2cm,right=2.2cm]{geometry}
\usepackage{amsfonts}
\usepackage{mathrsfs,amscd,amssymb,amsthm,amsmath,bm,graphicx,psfrag,subfigure,url,mathtools}
\usepackage{pict2e}
\usepackage{stfloats}
\usepackage{psfrag,amsmath}
\usepackage{tikz}
\usetikzlibrary{calc}
\usepackage{indentfirst}
\usepackage{hyperref}
\usepackage{bookmark}
\usepackage{enumerate}
\usepackage{latexsym,euscript,epic,eepic,color}
\usepackage{multirow}
\usepackage{multicol}
\usepackage{longtable}
\usepackage{adjustbox}
\usepackage[all]{xy}
\usepackage{setspace}
\usepackage{epstopdf}
\allowdisplaybreaks
\usepackage{authblk}
\usepackage{cmap}
\usepackage{pifont}
\usepackage{epstopdf}
\usepackage{caption}
\usepackage{makecell}
\usepackage{array}
\usepackage{bm}
\usepackage{psfrag}
\usepackage{verbatim}
\usepackage{xcolor}

\makeatletter

\renewcommand{\@seccntformat}[1]{{\csname the#1\endcsname}{\normalsize .}\hspace{.5em}}
\makeatother

\usepackage{ifpdf}
\usepackage{indentfirst}

\def \[{\begin{equation}}
\def \]{\end{equation}}

\usepackage{epstopdf}
\numberwithin{equation}{section}
\newtheorem{thm}{Theorem}[section]

\newtheorem{lem}[thm]{Lemma}

\newtheorem{remark}{Remark}[section]

\allowdisplaybreaks

\begin{document}
\captionsetup[figure]{labelfont={bf},name={Fig.},labelsep=period}
\baselineskip=0.23in

\title{\bf  Characterizing $A_\alpha$-minimizer graphs: given order and independence number}

\author[1]{Jiaqi Zhang}
\author[,1,2]{Shuchao Li\thanks{{\it Email address}: lscmath@ccnu.edu.cn (S. Li)}}
\affil[1]{School of Mathematics and Statistics, and Hubei Key Lab--Math. Sci.,\linebreak Central China Normal University, Wuhan 430079, China}
\affil[2]{Key Laboratory of Nonlinear Analysis \& Applications (Ministry of Education),\linebreak Central China Normal University, Wuhan 430079, China}
\date{\today}
\maketitle

\begin{abstract}
For a given graph \( G \), let \( A(G) \), \( Q(G) \), and \( D(G) \) denote the adjacency matrix, signless Laplacian matrix, and diagonal degree matrix of \( G \), respectively. The \( A_\alpha(G) \) matrix, proposed by Nikiforov, is defined as \( A_\alpha(G)=\alpha D(G)+(1 - \alpha)A(G) \), where \( \alpha\in[0,1] \). This matrix captures the gradual transition from \( A(G) \) to \( Q(G) \). Let \( \mathcal{G}_{n,\gamma} \) denote the family of all connected graphs with \( n \) vertices and independence number \( \gamma \). A graph in \( \mathcal{G}_{n,\gamma} \) is referred to as an \( A_\alpha \)-minimizer graph if it achieves the minimum \( A_\alpha \) spectral radius. In this paper, we first demonstrate that the \( A_\alpha \)-minimizer graph in \( \mathcal{G}_{n,\gamma} \) must be a tree when \( \gamma\geq\left\lceil\frac{n}{2}\right\rceil \), and we provide several characterizations of such \( A_\alpha \)-minimizer graphs. We then specifically characterize the \( A_\alpha \)-minimizer graphs for the case \( \gamma = \left\lceil\frac{n}{2}\right\rceil + 1 \) when $n\geq 9$. Furthermore, we obtain a structural characterization for the \( A_\alpha \)-minimizer graph when \( \gamma=n - c \), where \( c\geq4 \) is an integer. 

\vskip 0.2cm
\noindent {\bf Keywords:} Independence number; $A_\alpha$-index;  $A_\alpha$-minimizer graph; Tree
\vspace{2mm}

\noindent {\bf AMS Subject Classification:} 05C75; 05C50
\end{abstract}

\setcounter{section}{0}
\section{\normalsize Introduction}\setcounter{equation}{0}
Let \(G = (V(G), E(G))\) be a graph, where \(V(G)\) and \(E(G)\) denote the vertex set and edge set, respectively. As usual, \(|G|\) or \(|V(G)|\) denotes the order of \(G\), and \(|E(G)|\) denotes its size. All graphs considered herein are connected, finite, and simple. Unless otherwise stated, we follow standard notation and terminology (see \cite{11}).  

For a vertex \(v \in V(G)\), let \(N_G(v)\) denote the set of neighbors of \(v\). The degree of \(v\) is \(d_G(v) := |N_G(v)|\). The maximum degree of \(G\) is \(\Delta(G) := \max\{d_G(u) \mid u \in V(G)\}\), and the minimum degree is \(\delta(G) := \min\{d_G(u) \mid u \in V(G)\}\).  

A path \(P = v_0v_1\cdots v_s\) (\(s \geq 1\)) is an internal path if \(d(v_0), d(v_s) \geq 3\) and \(d(v_i) = 2\) for \(i \in [1, s-1]\); every edge in an internal path is an internal edge. A vertex subset \(R\subseteq V(G)\) is an independent set if the induced subgraph on \(R\) is empty. The independence number \(\gamma(G)\) is the size of the largest independent set in \(G\).  

A vertex \(v\) in a tree \(T\) is a branch point if \(d_T(v) \geq 3\). An end branch point is a branch point that does not lie on any path connecting two other branch points.  

For a graph \(G\) of order \(n\), the adjacency matrix \(A(G) = (a_{ij})\) is an \(n \times n\) matrix where \(a_{ij} = 1\) if \(v_i\) is adjacent to \(v_j\), and \(0\) otherwise. The degree matrix \(D(G) = \text{diag}(d(v_1), d(v_2), \ldots, d(v_n))\) is diagonal. The signless Laplacian matrix is defined as \(Q(G) = D(G) + A(G)\).  

Nikiforov \cite{1} defined the \(A_{\alpha}\)-matrix in 2017 as:  
$$
A_{\alpha}(G) = \alpha D(G) + (1-\alpha) A(G), \quad 0 \leq \alpha \leq 1.
$$  
Clearly, \(A_0(G)\) is the adjacency matrix \(A(G)\), \(2A_{1/2}(G)\) equals the signless Laplacian \(Q(G)\), and \(A_1(G)\) is the degree matrix \(D(G)\). The \(A_{\alpha}(G)\) is real, symmetric, and non-negative, so its eigenvalues are real. Let \(\lambda_{\alpha}(G)\) denote its largest eigenvalue (the \(A_{\alpha}\) spectral radius or \(A_{\alpha}\) index). By the Perron-Frobenius theory, a connected graph \(G\) has a unit positive eigenvector \(\boldsymbol{x} = (x_1, \ldots, x_n)^T\) associated with \(\lambda_{\alpha}(G)\), called the Perron vector of \(A_{\alpha}(G)\), where \(x_i\) is the ``weight" of vertex \(v_i\).  

Nikiforov \cite{1} introduced \(A_{\alpha}(G)\) to track the transition from \(A(G)\) to \(Q(G)\), and this matrix has since attracted significant attention (see \cite{LL2025,4,YGL2024,ZL2025} and references therein).  

Let \(\mathcal{G}_{n,\gamma}\) be the family of connected \(n\)-vertex graphs with independence number \(\gamma\). An \(A_{\alpha}\)-minimizer graph in \(\mathcal{G}_{n,\gamma}\) minimizes \(\lambda_{\alpha}(G)\). Research on \(A_0\)- and \(A_{1/2}\)-minimizer graphs is extensive: Xu et al. \cite{9} characterized graphs in \(\mathcal{G}_{n,\gamma}\) with the smallest adjacency spectral radius for \(\gamma \in \{n-1, n-2, n-3, \lceil \frac{n}{2} \rceil + 1, \lceil \frac{n}{2} \rceil, 2, 1\}\); Lou and Guo \cite{7} showed that \(A_0\)-minimizers in \(\mathcal{G}_{n,\gamma}\) are trees if \(\gamma \geq \lceil \frac{n}{2} \rceil\); Hu et al. \cite{8} refined these results for the adjacency spectral radius; and Hu et al. \cite{12} proved analogous results for \(A_{1/2}\)-minimizers.  

Sun et al. \cite{4} characterized \(A_{\alpha}\)-minimizers in \(\mathcal{G}_{n,\gamma}\) for \(\alpha \in [0,1)\) and \(\gamma \in \{1, \lfloor \frac{n}{2} \rfloor, \lfloor \frac{n}{2} \rfloor + 1, n-3, n-2, n-1\}\). Building on [\cite{8},\cite{12},\cite{7},\cite{4},\cite{9}], we focus on \(A_{\alpha}\)-minimizers in \(\mathcal{G}_{n,\gamma}\).  

\begin{figure}[htbp]
\centering
\includegraphics[scale=0.90]{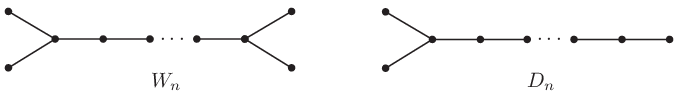}
\caption{Graphs $W_n$ and $D_n$.}
\label{Fig.1}
\end{figure}

Sun et al. \cite{4} handled \(\gamma \in \{n-2, n-3\}\) for \(A_{\alpha}\)-minimizers, so we naturally extend to \(\gamma \geq \lceil \frac{n}{2} \rceil\). Our first main result is:
\begin{thm}\label{thm3.01}
Let \(\gamma \geq \lceil \frac{n}{2} \rceil\) and \(\alpha \in [0,1)\). If \(G^*\) is an \(A_{\alpha}\)-minimizer in \(\mathcal{G}_{n,\gamma}\), then \(G^*\) is a tree.\end{thm}
Our second result specifies \(A_{\alpha}\)-minimizers for \(\gamma = \lceil \frac{n}{2} \rceil + 1\):  
\begin{thm}\label{thm3.02}
Let \(H\) be an \(A_{\alpha}\)-minimizer in \(\mathcal{G}_{n, \lceil \frac{n}{2} \rceil + 1}\) with \(\alpha \in [0,1)\) and \(n \geq 9\). Then \(H \cong D_n\) if \(n\) is even, and \(H \cong W_n\) if \(n\) is odd (see \(D_n\) and \(W_n\) in Fig.~\ref{Fig.1}).
\end{thm}
The subdivision graph \(S(T)\) of a tree \(T\) is formed by subdividing every edge of \(T\) once. 
Our third result characterizes \(A_{\alpha}\)-minimizers for \(\gamma \in [\lceil \frac{n}{2} \rceil + 2, n-4]\). 
\begin{thm}\label{thm4.1}
Let \(c \geq 4\) be an integer, \(\gamma = n - c \in [\lceil \frac{n}{2} \rceil + 2, n-4]\), and \(\hat{T}\) be an \(A_{\alpha}\)-minimizer in \(\mathcal{G}_{n,\gamma}\) with \(0 \leq \alpha < 1\). Let \(T\) be a \((n-\gamma)\)-vertex tree, and \(S(T)\) its subdivision graph. Then \(\hat{T}\) is obtained by attaching \(2\gamma - n + 1\) pendant edges to vertices in \(V(T)\) if one of the following holds:  
\begin{enumerate}[{\rm (i)}]
    \item \(\alpha = 0\) and \(n \geq \lceil \frac{4c^2 - 6c - 3}{3} \rceil;\)
    \item \(\alpha \in (0, \frac{1}{2})\) and \(n \geq \max\left\{ c\left( \frac{(1-2\alpha)(c-3)}{3\alpha} + \sqrt{c+2} \right)^2 - 3c - 1, 2c + 4 \right\};\)  
    \item \(\alpha \in [\frac{1}{2}, 1)\) and \(n \geq \max\{ c^2 - c - 1, 2c + 4 \}\). 
\end{enumerate}
\end{thm}

\begin{figure}[htbp]
\centering
\includegraphics[scale=0.9]{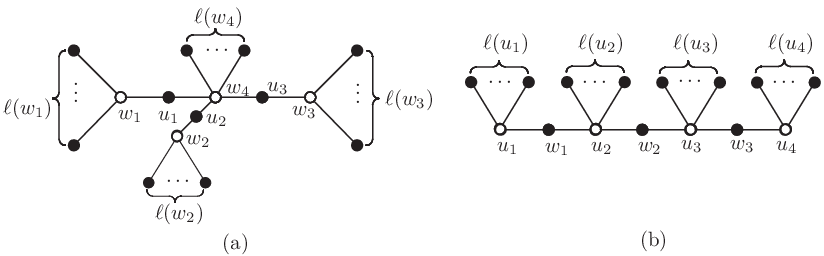}
\caption{(a) $S(K_{1,3})\circ_{V(K_{1,3})}(\ell(w_1),\ell(w_2),\ell(w_3),\ell(w_4))$; (b) $S(P_4)\circ_{V(P_4)}(\ell(u_1),\ell(u_2),\ell(u_3),\ell(u_4))$.}
\label{Fig.4}
\end{figure}

All the possible structures of the \(A_{\alpha}\)-minimizer graphs in \(\mathcal{G}_{n,n-4}\) under the conditions of Theorem \ref{thm4.1} are depicted in Fig.~\ref{Fig.4}.

In the remainder of this paper, we present some preliminaries in Section \ref{s2}. In Section \ref{s3}, we give the proofs of Theorems~\ref{thm3.01} and \ref{thm3.02}. Moreover, we will give some further characterizations on the structure for all the $A_\alpha$-minimizer graphs with $\gamma\in[\lceil\frac{n}{2}\rceil+2,n-2]$ and $\alpha\in[0,1)$. Theorem~\ref{thm4.1} is proved in Section \ref{s4}. In the last section, leveraging the results obtained in this paper, we determine $18$ possible $A_\alpha$-minimizer graphs in $\mathcal{G}_{n,n-4}$ for $\frac{1}{2}\leq\alpha< 1$. This work refines and extends the findings of Liu and Wang \cite{13}.

\section{\normalsize Some known lemmas}\label{s2}\setcounter{equation}{0}
Some necessary preliminaries are given in this section, which we will use to prove our main results.

\begin{lem}[\cite{1}]\label{lem2.1}
For a connected graph \(G\) and \(\alpha \in [0,1)\), if \(H\) is a proper subgraph of \(G\), then \(\lambda_{\alpha}(H) < \lambda_{\alpha}(G)\). 
\end{lem}

\begin{lem}[\cite{2}]\label{lem2.2}
Let \(G\) be connected, \(\boldsymbol{x}\) the Perron vector of \(A_{\alpha}(G)\) (\(\alpha \in [0,1)\)), and \(u, v \in V(G)\) with \(R \subseteq N(v) \setminus (N(u) \cup \{u\})\). Let \(G^* = G - \{vw \mid w \in R\} + \{uw \mid w \in R\}\). If \(R \neq \emptyset\) and \(x_u \geq x_v\), then \(\lambda_{\alpha}(G) < \lambda_{\alpha}(G^*)\).  
\end{lem}


\begin{lem}[\cite{10}]\label{lem2.4}
For a connected graph \(G\), let \(G_{s,t}\) be formed by attaching pendant paths of lengths \(s, t\) to a vertex of \(G\). If \(1 \leq t \leq s\), then \(\lambda_{\alpha}(G_{s+1,t-1}) < \lambda_{\alpha}(G_{s,t})\) for \(0 \leq \alpha < 1\). 
\end{lem}

\begin{lem}[\cite{4}]\label{lem2.5}
For a connected \(n\)-vertex graph \(G\) (\(n \geq 9\)), if \(G \not\cong P_n\) or \(C_n\), then \(\lambda_{\alpha}(D_n) \leq \lambda_{\alpha}(G)\) for \(0 \leq \alpha < 1\), with equality if and only if \(G \cong D_n\) (see\ Fig.~\ref{Fig.1}).
\end{lem}

Let $T$ be a tree, denote by $L(T)$ the set of all leaves of $T$.

\begin{lem}[\cite{5}]\label{lem2.6}
Every tree \(T\) has a maximum independent set \(I(T)\) containing all leaves of \(T\).
\end{lem}


\begin{lem}[\cite{12}]\label{lem3.7}
For a tree \(T\) with independence number \(\gamma(T)\) and a maximum independent set \(I(T)\) containing all leaves, if vertices on the path between any two leaves alternate in \(I(T)\), then \(\gamma(T) \geq \frac{n+1}{2}\).  
\end{lem}

\begin{lem}[\cite{6}]\label{lem2.7}
Let \(uv\) be an internal edge of \(G\). If \(G'\) is formed by subdividing \(uv\), then \(\lambda_{\alpha}(G') \leq \lambda_{\alpha}(G)\) for \(0 \leq \alpha < 1\).  
\end{lem}

\begin{lem}[\cite{7}]\label{lem2.8}
Let \(uv\) be an edge of \(G\). Let \(G_1\) (resp. \(G_2\)) be formed by subdividing \(uv\) once (resp. twice). Then: 
 \begin{enumerate}[{\rm(i)}]
      \item $\gamma(G_1)\in\{\gamma(G),\gamma(G)+1\}$,
      \item $\gamma(G_2)=\gamma(G)+1$.
    \end{enumerate}
\end{lem}

For a subset \(Y = \{u_1, \ldots, u_t\} \subseteq V(G)\), let \(G \circ_Y(\ell_1, \ldots, \ell_t)\) denote the graph formed by attaching \(\ell_i\) pendant edges to \(u_i\) (\(i = 1, \ldots, t\)); \(Y\) is the rooted set. We write \(G \circ_Y(k)\) if \(\ell_1 = \cdots = \ell_t = k\).  

\begin{lem}[\cite{16}]\label{lem4.0}
Let \(A\) be a partite set of a bipartite graph \(G\), and \(G' = G \circ_A(k)\). Then \(\lambda_0(G') = \sqrt{\lambda_0^2(G) + k}\). 
\end{lem}

\begin{lem}[\cite{3}]\label{lem4.1}
For a graph \(G\) and \(\alpha \in [0, \frac{1}{2}]\),  
$
\lambda_{\alpha}(G) \leq 2\alpha \lambda_{1/2}(G) + (1-2\alpha)\lambda_0(G).
$
Equality holds for connected irregular \(G\) if and only if \(\alpha = 0\) or \(\alpha = \frac{1}{2}\).  
\end{lem}

\begin{lem}[\cite{17}]\label{lem4.2.2}
The path \(P_n\) uniquely minimizes the \(A_0\) spectral radius among \(n\)-vertex connected graphs, with \(\lambda_0(P_n) = 2\cos\left( \frac{\pi}{n+1} \right)\).  \end{lem}

\begin{lem}[\cite{1}]\label{lem4.2}
For a graph \(G\) with \(\Delta(G) = \Delta\) and \(\alpha \in [0,1)\),  
$
\lambda_{\alpha}(G) \geq \frac{1}{2}(\sqrt{\alpha^2(\Delta + 1)^2 + 4\Delta(1-2\alpha)}\linebreak + \alpha(\Delta + 1)).
$
Equality holds for connected \(G\) if and only if \(G = K_{1,\Delta}\). In particular,
$$
\lambda_\alpha(G) \ge 
\begin
{cases}
\alpha(\Delta + 1), & \text{if } \alpha \in [0, \frac{1}{2}]; \\
\alpha\Delta + \frac{(1 - \alpha)^2}{\alpha}, & \text{if } \alpha \in [\frac
{1}{2}, 1).
\end
{cases}
$$
\end{lem}

\begin{lem}[\cite{1}]\label{lem4.2.1}
For \(0 \leq \alpha < 1\),  
$$
\frac{2|E(G)|}{|V(G)|} \leq \lambda_{\alpha}(G) \leq \max_{uv \in E(G)} \{ \alpha d_G(u) + (1-\alpha)d_G(v) \}.
$$  
Left equality holds if and only if \(G\) is regular.  
\end{lem}

Recall that the subdivision graph \(S(T)\) of a tree \(T\) is formed by subdividing every edge of \(T\) once.
\begin{lem}[\cite{12}]\label{lem4.3}
A tree \(G\) is isomorphic to \(S(T)\) (for some tree \(T\) on \(t\) vertices) iff \(G\) is bipartite with partition \((R, S)\) where all vertices in \(R\) have degree 2 and \(|S| = t\).  
\end{lem}

\begin{lem}[\cite{4}]\label{lem3.1}
For an \(A_{\alpha}\)-minimizer \(G\) in \(\mathcal{G}_{n,\gamma}\): if \(\gamma(G) = \lceil \frac{n}{2} \rceil\), then \(G \cong P_n\); if \(\gamma(G) = n-1\), then \(G \cong K_{1,n-1}\). 
\end{lem}

\begin{lem}[\cite{4}]\label{lem5.1}
Let \(u, v \in V(G)\) with \(G - u \cong G - v\). Let \(G(s,t)\) be formed by attaching \(s\) (resp. \(t\)) pendant edges to \(u\) (resp. \(v\)) (\(s \geq t \geq 1\)). For \(\alpha \in [0,1)\), \(\lambda_{\alpha}(G(s,t)) < \lambda_{\alpha}(G(s+1,t-1))\).  
\end{lem}
\begin{lem}[\cite{13}]\label{lem5.2}
Let \(G = S(P_4) \circ_{V(P_4)}(\ell(u_1), \ldots, \ell(u_4))\) (Fig.~\ref{Fig.4}(b)) with \(\max_{1 \leq i,j \leq 4} |d_G(u_i) - d_G(u_j)| \leq 2\). Let \(x\) be the Perron vector of \(A_{\alpha}(G)\) for \(\frac{1}{2} \leq \alpha < 1\). If \(\ell(u_1) \geq 1\) (resp. \(\ell(u_4) \geq 1\)) and \(\Delta(G) > \ell(u_1) + 1\) (resp. \(\Delta(G) > \ell(u_4) + 1\)), then \(x_{u_1} \leq x_{u_2}\) (resp. \(x_{u_4} \leq x_{u_3}\)). 
\end{lem}
For a matrix \(W\) indexed by \(S = \{1, \ldots, n\}\) and a partition \(\pi = S_1 \cup \cdots \cup S_t\) of \(S\), the quotient matrix \(W_{\pi}\) is a \(t \times t\) matrix with entries equal to the average row sums of submatrices \(W_{ij}\) (indexed by \(S_i, S_j\)). If row sums of \(W_{ij}\) are constant, \(\pi\) is equitable, and \(W_{\pi}\) is the equitable quotient matrix.  
\begin{lem}[\cite{14}]\label{lem5.3}
For an equitable partition \(\pi\) of a non-negative matrix \(W\), the largest eigenvalue of \(W\) equals that of its equitable quotient matrix \(W_{\pi}\).  
\end{lem}

\begin{lem}[\cite{13}]\label{lem5.4}
Let \(T_1 \cong S(P_4) \circ_{V(P_4)}(t,t,t,t)\) and \(T_2 \cong S(P_4) \circ_{V(P_4)}(t+1,t-1,t-1,t+1)\). For \(\alpha \in [\frac{1}{2}, 1)\), \(\lambda_{\alpha}(T_1) > \lambda_{\alpha}(T_2)\).  
\end{lem}

\begin{lem}[\cite{9}]\label{lem5.9.1}
For an \(n\)-vertex connected graph \(G\) with \(\gamma(G) = 2\), \(|E(G^c)| \leq \lfloor \frac{n^2}{4} \rfloor - 1\).
\end{lem}
Vertices \(u, v \in V(G)\) are \textit{equivalent} if there is an automorphism \(\varphi: G \to G\) with \(\varphi(u) = v\).
\begin{lem}[\cite{1}]\label{lem5.9.2}
For a connected \(n\)-vertex graph \(G\) and Perron vector \(\boldsymbol{x}\) of \(A_{\alpha}(G)\), equivalent vertices have equal weights in \(\boldsymbol{x}\).
\end{lem}

\section{\normalsize Proofs of Theorems~\ref{thm3.01} and \ref{thm3.02}}\label{s3}\setcounter{equation}{0}
Theorems~\ref{thm3.01} and \ref{thm3.02} are proved in this section. Theorem~\ref{thm3.01} shows that the graph among $\mathcal{G}_{n,\gamma}$ having minimal $A_\alpha$ spectral radius is a tree for  $\gamma\geq\lceil\frac{n}{2}\rceil$ and $\alpha\in[0,1)$, and Theorem~\ref{thm3.02} specifically characterizes the $A_\alpha$-minimizer graphs among $\mathcal{G}_{n,\lceil\frac{n}{2}\rceil+1}$ with given order $n\geq 9$ and $\alpha\in[0,1)$.

\begin{proof}[\bf Proof of Theorem \ref{thm3.01}]
By Lemma \ref{lem3.1}, the result holds for \(\gamma \in \{\lceil \frac{n}{2} \rceil, n-1\}\). Consider \(G^* \in \mathcal{G}_{n,\gamma}\) with \(\lceil \frac{n}{2} \rceil < \gamma(G^*) < n-1\). If \(G^*\) has cycles, we can find a spanning tree, say \(\hat{T}\) other than $K_{1,n-1}$ of \(G^*\). By Lemma \ref{lem2.1}, \(\lambda_{\alpha}(\hat{T}) < \lambda_{\alpha}(G^*)\). Obvoiusly, any two non-adjacent vertices of $G^*$ are still not adjacent in the spanning tree \(\hat{T}\) of \(G^*\), thus \(\gamma(\hat{T}) \geq \gamma(G^*)\). If \(\gamma(\hat{T}) = \gamma(G^*)\), \(G^*\) is not minimal, so:  
\begin{align}\label{3.1}
\gamma(\hat{T})>\gamma(G^*)>\lceil\frac{n}{2}\rceil ~\mbox{and}~\lambda_\alpha(\hat{T})<\lambda_\alpha(G^*).
\end{align} 
\(\hat{T} \not\cong P_n\) (else \(\gamma(\hat{T}) = \lceil \frac{n}{2} \rceil\), contradicting \eqref{3.1}). We analyze two cases: 

\begin{figure}
\centering
\includegraphics[scale=0.86]{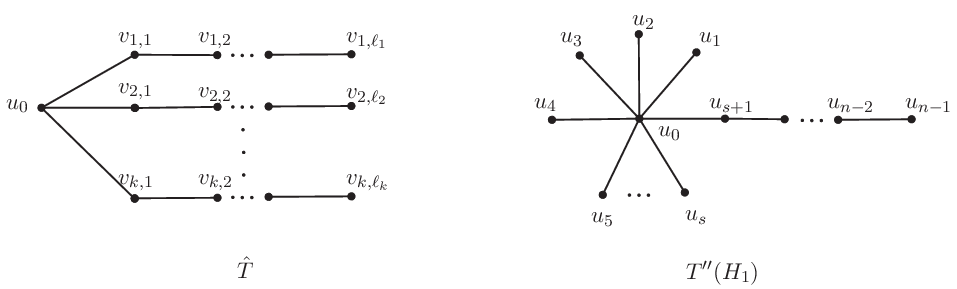}
\caption{Two graphs $\hat{T},\,T''(H_1)$ in the proofs of Theorems \ref{thm3.01} and \ref{thm3.02}.}
\label{Fig.2}
\end{figure}

{\bf Case 1.} $\hat{T}$ contains a unique branch point, say $u_0$. Let $\hat{T}$ be a tree obtained by attaching $k$ paths $P_{{\ell_1}+1}, P_{{\ell_2}+1},\dots, P_{{\ell_k}+1}$ to vertex $u_0$, where $k\geq3$ and $\ell_i\geq\ell_j$ for $1\leq i\leq j\leq k$. For simplicity, let  $P_{{\ell_j}+1}=u_0v_{j,1}v_{j,2}\cdots v_{j,\ell_j}$ for $j=1,\dots, k$ (see Fig. \ref{Fig.2}). Note that $\hat{T}\not\cong K_{1,n-1}$, then we have $\ell_1\geq2$. Let
\[
T'=\begin{cases}
    \hat{T}&\mbox{if}\quad\ell_2=1,\\
    \hat{T}-u_0v_{2,1}+v_{1,\ell_1}v_{2,1}&\mbox{if}\quad\ell_2\geq2\quad\mbox{is even},\\
    \hat{T}-v_{2,1}v_{2,2}+v_{1,\ell_1}v_{2,2}&\mbox{if}\quad\ell_2\geq3\quad\mbox{is odd}.
\nonumber
\end{cases}
\]
Clearly, $|T'|=|\hat{T}|$. According to Lemma \ref{lem2.4}, one sees $\lambda_\alpha(T')\leq\lambda_\alpha(\hat{T})$. By Lemmas \ref{lem2.6} and \ref{lem2.8}, we get  $\gamma(T')=\gamma(\hat{T})$. Repeat the above process finitely many times on $T'$, 
we ultimately get a tree $T''$ so that
\begin{align}\label{3.2}
|T''|=|T'|=|\hat{T}|,
\gamma(T'')=\gamma(T')=\gamma(\hat{T})
~\mbox{and}~
\lambda_\alpha(T'')\leq\lambda_\alpha(T')\leq\lambda_\alpha(\hat{T}),
\end{align}
where $T''$ is constructed by attaching a pendant path to the center vertex of $K_{1,s}$  (see Fig. \ref{Fig.2}). 

When $s=1$, one gets $T''\cong P_n$. By \eqref{3.2} we have $\gamma(\hat{T})=\gamma(P_n)=\lceil\frac{n}{2}\rceil$, a contradiction to \eqref{3.1}. 

When $s=2$, one gets $T''\cong D_n$ (see Fig. \ref{Fig.1}). By \eqref{3.1} and \eqref{3.2} one gets 
\[\label{3.03}
 \gamma(\hat{T})=\gamma(D_n)=\lceil\frac{n+1}{2}\rceil>\gamma(G^*)>\lceil\frac{n}{2}\rceil.
\]
Note that $\lceil\frac{n+1}{2}\rceil=\lceil\frac{n}{2}\rceil$ holds for odd $n$, and $\lceil\frac{n+1}{2}\rceil=\lceil\frac{n}{2}\rceil+1$ holds for even $n$. This gives an obvious contradiction to \eqref{3.03}.

If $s=3$, then again by \eqref{3.1} and \eqref{3.2} one sees 
$$\gamma(\hat{T})=\gamma(T'')=\lceil\frac{n}{2}\rceil+1>\gamma(G^*)>\lceil\frac{n}{2}\rceil,$$
also a contradiction.


So in what follows, we consider $s\geq 4$. In this case, let $T'''=T''-u_0u_1+u_{n-2}u_1$. Then $|T''|=|T'''|$ and $\gamma(T''')=\gamma(T'')$. Assume ${\bf x}=(x_{u_0},\dots, x_{u_{n-1}})^T$ is the $A_\alpha(T''')$'s Perron vector.

If $x_{u_{n-2}}\leqslant x_{u_0}$, then by Lemma \ref{lem2.2}, one has 
$\lambda_\alpha(T''')<\lambda_\alpha(T'''-u_{n-2}u_1+u_0u_1)=\lambda_\alpha(T'').$
If $x_{u_{n-2}}> x_{u_0}$, then due to Lemma \ref{lem2.2}, one also has 
$\lambda_\alpha(T''')<\lambda_\alpha(T'''-\{u_0u_i|i=2,\dots,s-1\}+\{u_{n-2}u_i|i=2,\dots,s-1\})=\lambda_\alpha(T'').$
Bearing in mind with \eqref{3.1} and \eqref{3.2} gives us
\begin{align}\label{3.3}
\lambda_\alpha(T''')<\lambda_\alpha(T'')\leq\lambda_\alpha(\hat{T})<\lambda_\alpha(G^*).
\end{align}

Let $T'''_1:=T'''$, and, for $2\leq i\leq s-2,$ let $T'''_i$ be constructed from $T'''_{i-1}$ by subdividing an internal edge and deleting a leaf vertex $u_i$ simultaneously. Clearly, for $i=2,\dots,s-2$, one gets $|T'''_i|=|T'''_{i-1}|$ and $T'''_{s-2}\cong W_n$ (see Fig. \ref{Fig.1}).

By Lemmas \ref{lem2.6} and \ref{lem2.8}, for $2\leq i\leq s-2$, we get  
$
\gamma(T'''_i)\in\{\gamma(T'''_{i-1}),\gamma(T'''_{i-1})-1\}. 
$ By Lemmas \ref{lem2.1} and \ref{lem2.7}, for $2\leq i\leq s-2$, we have 
\begin{align}\label{3.4}
\lambda_\alpha(T'''_i)<\lambda_\alpha(T'''_{i-1}).
\end{align}

If $\gamma(T'''_j)>\gamma(G^*)$ for all $j\in\{1,\dots,s-2\}$,
then by \eqref{3.1}, we obtain
${\lceil\frac{n}{2}\rceil}+1=\gamma(W_n)=\gamma(T'''_{s-2})\geq\gamma(G^*)+1>{\lceil\frac{n}{2}\rceil}+1,$
which is impossible. Hence, there is an $j\in\{1,2,\dots,s-2\}$ so that $\gamma(T'''_j)=\gamma(G^*)$.
Combining with \eqref{3.3} and \eqref{3.4}, we have
$$\lambda_\alpha(T'''_j)\leq\lambda_\alpha(T'''_1)<\lambda_\alpha(T'')\leq\lambda_\alpha(\hat{T})<\lambda_\alpha(G^*),$$
also a contradiction to the minimality of $\lambda_\alpha(G^*)$.

{\bf Case 2.} $\hat{T}$ contains at least 2 branch points. Let $u,v$ be two branch points of 
$\hat{T}$ such that the unique path $P$ between $u$ and $v$ is an internal path. Since $u,v$ are branch points, we have $d_{\hat{T}}(u),d_{\hat{T}}(v)\geq 3$. Let $\{u_1,u_2\}\subseteq N_{\hat{T}}(u)\setminus V(P), \{v_1,v_2\}\subseteq N_{\hat{T}}(v)\setminus V(P)$.

Let $T^\diamond_0=\hat{T}$ and, for $i=0,1,2,\ldots, k-2,$ subdividing an internal edge on the internal path between $u,v$ and deleting a leaf vertex $w$ of $T^\diamond_{i}$ ($w\not \in \{u_1, u_2, v_1, v_2\}$) simultaneously yields the graph $T^\diamond_{i+1};$ see Fig.~\ref{F8}. 
\begin{figure}
  \centering
  \includegraphics[width=125mm]{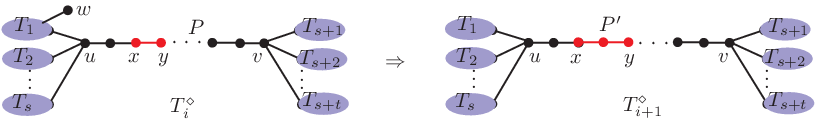}
  \caption{Graphs $T^\diamond_i$ and $T^\diamond_{i+1}$}\label{F8}
\end{figure}

Repeating this process yields a tree sequence $T^\diamond_0,T^\diamond_1,T^\diamond_2,\dots,T^\diamond_{k-1}$ such that $T^\diamond_{k-1}\cong W_n$. Subdividing an internal edge of $T^\diamond_{k-1}$ and deleting one of its leaf vertex in $\{u_1, u_2, v_1, v_2\}$ yields $T^\diamond_{k}$. Clearly, we obtain $T^\diamond_{k}\cong D_n$.

By Lemmas \ref{lem2.1} and \ref{lem2.7}, for $0\leqslant i\leqslant k-1$, we have 
$
\lambda_\alpha(T^\diamond_{i+1})<\lambda_\alpha(T^\diamond_i)<\lambda_\alpha(G^*).
$ By Lemmas \ref{lem2.6} and \ref{lem2.8}, for $0\leqslant i\leqslant k-1$, we have
\begin{align*}
\gamma(T^\diamond_0)=\gamma(\hat{T})>\gamma(G^*) ~\mbox{and}~ \gamma(T^\diamond_{i})-1\leq\gamma(T^\diamond_{i+1})\leq\gamma(T^\diamond_{i})+1.
\end{align*}

If $\gamma(T^\diamond_j)=\gamma(G^*)$ holds for $j\in [1,k]$, then it contradicts the minimality of $G^*$. 
So we have $\gamma(T^\diamond_j)>\gamma(G^*)$ for all $j\in\{1,2,\dots,k\}.$ Particularly, $\gamma(T^\diamond_k)>\gamma(G^*)>\lceil\frac{n}{2}\rceil$, a contradiction to $\gamma(T^\diamond_k)=\gamma(D_n)=\lceil\frac{n+1}{2}\rceil$.
\end{proof}
\begin{proof}[\bf Proof of Theorem \ref{thm3.02}]
For even $n$, Lemma \ref{lem2.5} shows that $D_n$ is the unique $A_\alpha$-minimizer graph with independence number  ${\lceil\frac{n+1}{2}\rceil}={\lceil\frac{n}{2}\rceil}+1$. Thus $H\cong D_n$. If $n$ is odd, then we can easily show that $H\not\cong P_n, K_{1,n-1},$ and $H$ is a tree (based on Theorem \ref{thm3.01}). Next, let's consider the remaining two cases.

{\bf Case 1.} $H$ contains a unique branch point.
By a similar discussion as Case 1 in the proof of Theorem \ref{thm3.01}, we obtain a tree being isomorphic to $H_1$ such that
\begin{align}\label{3.5}
\ |H_1|=|H|,
\gamma(H_1)=\gamma(H)
~\mbox{and}~
\lambda_\alpha(H_1)\leq\lambda_\alpha(H).
\end{align}
Recall that $H_1$ is a tree which is formed from $K_{1,s}$ by means of attaching a pendant path to its center vertex (see Fig. \ref{Fig.2}). 

If $s=1$, then $H_1\cong P_n$. By \eqref{3.5}, we have $\gamma(H)=\gamma(H_1)=\lceil\frac{n}{2}\rceil$, a contradiction to the fact that $\gamma(H)=\lceil\frac{n}{2}\rceil+1$. 

If $s=2$, then $H_1\cong D_n$ (see Fig. \ref{Fig.1}). By \eqref{3.5}, we have 
$\gamma(H)=\gamma(H_1)={\lceil\frac{n+1}{2}\rceil}$.
Note that $\lceil\frac{n+1}{2}\rceil=\lceil\frac{n}{2}\rceil$ if $n$ is odd. This also gives an obvious contradiction to the definition of $H$.

So in what follows, we consider $s\geq3$. In this case, let $H_2=H_1-u_0u_1+u_{n-2}u_1$. Then $|H_2|=|H_1|$ and $\gamma(H_2)=\gamma(H_1)$. 
The $A_\alpha(H_2)$'s Perron vector is denoted by  ${\bf x}=(x_{u_0},\dots, x_{u_{n-1}})^T$.

By virtue of Lemma \ref{lem2.2}, if $x_{u_{n-2}}\leq x_{u_0}$, then we have 
$\lambda_\alpha(H_2)<\lambda_\alpha(H_2-u_{n-2}u_1+u_0u_1)=\lambda_\alpha(H_1).$ Similarly, 
if $x_{u_{n-2}}> x_{u_0}$, then $\lambda_\alpha(H_2)<\lambda_\alpha(H_2-\{u_0u_i|i=2,3,\dots,s-1\}+\{u_{n-2}u_i|i=2,3,\dots,s-1\})=\lambda_\alpha(H_1).$
Together with \eqref{3.5}, we obtain $\lambda_\alpha(H_2)<\lambda_\alpha(H_1)\leq\lambda_\alpha(H)$ and $\gamma(H_2)=\gamma(H_1)=\gamma(H)$,
also a contradiction to the minimality of $\lambda_\alpha(H)$.

{\bf Case 2.} $H$ contains at least 2 branch points. Denote with $H^\diamond_0=H$. With a discussion similar to Case 2 in Theorem \ref{thm3.01}’s proof, there exists a tree sequence $H^\diamond_0,H^\diamond_1,H^\diamond_2,\dots,H^\diamond_{k-1},H^\diamond_k$ such that $H^\diamond_{k-1}\cong W_n$ and $H^\diamond_{k}\cong D_n$.

By Lemmas \ref{lem2.1} and \ref{lem2.7}, for $0\leq i\leq k-1$, one has $\lambda_\alpha(H^\diamond_{i+1})<\lambda_\alpha(H^\diamond_i).$ 
By virtue of Lemmas \ref{lem2.6} and \ref{lem2.8}, one sees $\gamma(H^\diamond_{i})-1\leq\gamma(H^\diamond_{i+1})\leq\gamma(H^\diamond_{i})+1$ for $0\leq i\leq k-1$.

If $H\not\cong W_n$, then $\gamma(W_n)=\gamma(H)={\lceil\frac{n}{2}\rceil}+1$ and $\lambda_\alpha(W_n)<\lambda_\alpha(H)$,
a contradiction to the minimality of $\lambda_\alpha(H)$. Hence $H\cong W_n$.

Together with Cases 1 and 2, we complete our proof.  
\end{proof}
\section{\normalsize Proof of Theorem~\ref{thm4.1}}\label{s4}\setcounter{equation}{0}
Within this section we'll explore the proof for Theorem~\ref{thm4.1}. In order to do so, we need some structure lemmas. Let $\hat{T}$ have the smallest $A_\alpha$ spectral radius within $\mathcal{G}_{n,\gamma}$, where $\alpha\in[0,1)$ and ${\lceil\frac{n}{2}\rceil}+2\leq\gamma\leq n-2$.
\subsection{\normalsize Some structure characterization on $\hat{T}$}
In this subsection we need the following lemmas to characterize some structure of $\hat{T}$, which will be used in the proof of Theorem~\ref{thm4.1}.
\begin{lem}\label{lem3.0}
For any tree $T$ with at least $3$ branch points, it contains at least an end branch point.
\end{lem}
\begin{proof}
Suppose that $T$ has no end branch point. Then for any branch point $u$ in $T$, there exist other two branch points $u_1, u_2$, such that $u$ lies on the path that joins $u_1$ and $u_2$ in $T$. This leads to a contradiction since $T$ is a tree and the number of the branch points of $T$ is finite.
\end{proof}

\begin{lem}\label{Prop3.1}
$\hat{T}$ contains at least $2$ end branch points.
\end{lem}
\begin{proof}
We firstly show that $\hat{T}$ contains at least 2 branch points. Otherwise, assume that $\hat{T}$ contains a unique branch point (see Fig. \ref{Fig.2}). By means of a discussion analogous to that of Case 1 in Theorem \ref{thm3.02}, there exists a tree $T'$ such that $|T'|=|\hat{T}|,\gamma(T')=\gamma(\hat{T})$ and $\lambda_\alpha(T')<\lambda_\alpha(\hat{T})$, which contradicts the minimality of $\lambda_\alpha(\hat{T})$. Hence, $\hat{T}$ contains at least 2 branch points.

Denote $\tau(\hat{T})$ the number of the branch points in $\hat{T}$. Then we establish by induction on $\tau(\hat{T})$, that $\hat{T}$ contains at least 2 end branch points. If $\tau(\hat{T})=2$, then these two branch points are exactly the two end branch points in $\hat{T}$. 
By Lemma~\ref{lem3.0}, one sees that $\hat{T}$ contains at least 1 end branch point. Choose an end branch point, say $u$, in $\hat{T}$ such that there exist exactly $d_{\hat{T}}(u)-1$ pendant paths attached to it.

By deleting $d_{\hat{T}}(u)-1$ pendant paths mentioned above, we can obtain a new graph $T_1$ from $\hat{T}$, where $d_{T_1}(u)=1$, $d_{T_1}(v)=d_{\hat{T}}(v)$ for $v\in V(T_1)\setminus\{u\}$ and $\tau(T_1)=k-1$. By induction, $T_1$ has at least two end branch points $u_1,u_2$. Consider in $T_1$, it is straightforward to check that $u_1$ attaches $k=d_{T_1}(u_1)-1$ pendant paths $P_{11}, P_{12},\dots,P_{1k}$ and $u_2$ attaches $s=d_{T_1}(u_2)-1$ pendant paths $P_{21}, P_{22},\dots,P_{2s}$ ($k,s\geq2$).
In the following, we characterize the local structure of $\hat{T}$.

If $u$ is not on the $k+s$ paths mentioned above in $T_1$, then $u_1$ and $u_2$ still attach $k=d_{\hat{T}}(u_1)-1$ and $s=d_{\hat{T}}(u_2)-1$ pendant paths in $\hat{T}$, respectively. Thus $u$, $u_1$ and $u_2$ are end branch points in $\hat{T}$. Otherwise, for simplicity, assume that $u$ is on one of $k$ pendant paths attached to vertex $u_1$ in $T_1$. Then there are  $s=d_{\hat{T}}(u_2)-1$ pendant paths attached to $u_2$ in $\hat{T}$. Thus $u$ and $u_2$  are end branch points in $\hat{T}$. Therefore, $\hat{T}$ contains at least 2 end branch points. 
\end{proof}

\begin{lem}\label{Prop3.2}
For any end branch point $u_0$ of $\hat{T}$, there are exactly $d_{\hat{T}}(u_0)-1$ leaves in $N_{\hat{T}}(u_0)$.
\end{lem}
\begin{proof}
In view of Lemma \ref{Prop3.1}, assume that $u_0$ is an end branch point of $\hat{T}$. Then there exist exactly $d_{\hat{T}}(u_0)-1$ pendant paths attached to $u_0$. Denote by $w$ the remaining vertex which is adjacent to $u_0$. We proceed by considering the length of these pendant paths. 

Suppose that  $P=u_0u_1\cdots u_t$ is a pendant path attached to $u_0$ with $t\geq2$. According to Lemma \ref{Prop3.1}, $\hat{T}$ contains at least 2 branch points. Obviously, $u_0w$ must belong to an internal path of $\hat{T}$, which implies that $u_0w$ serves as an internal edge of $\hat{T}$. Let $T_1$ be the tree obtained from $\hat{T}$ by substituting the edge $u_0w$ by a path $u_0w_1w_2w$. Let $T_2=T_1-u_{t}-u_{t-1}$. By virtue of Lemmas \ref{lem2.1} and \ref{lem2.7}, one sees  
$\lambda_\alpha(\hat{T})\geq\lambda_\alpha(T_1)>\lambda_\alpha(T_2),~ |\hat{T}|=|T_1|-2=|T_2|.$
Moreover, according to Lemmas \ref{lem2.6} and \ref{lem2.8}, one gets  
$\gamma(\hat{T})=\gamma(T_1)-1=\gamma(T_2).$
Therefore, we obtain $T_2\in\mathcal{G}_{n,\gamma}$ and $\lambda_\alpha(T_2)<\lambda_\alpha(\hat{T})$, which leads contradiction on the minimality of $\lambda_\alpha(\hat{T})$. Then there are exactly $d_{\hat{T}}(u_0)-1$ leaves in $N_{\hat{T}}(u_0)$.
\end{proof}

\begin{lem}\label{Prop3.3}
Let $u_0,u_k$ be two end branch points of $\hat{T},$ and $P=u_0u_1\cdots u_k\,(k\geq1)$ be a path joining $u_0$ and $u_k$. Then $k$ is even and $u_i\in I(\hat{T})$ for odd $i\in[1,k-1]$.
\end{lem}
\begin{proof}
According to Lemma \ref{Prop3.2}, for any end branch point $u$ in $\hat{T}$, it attaches exactly $d_{\hat{T}}(u)-1$ leaf vertices. 
Let $v_0\in N_{\hat{T}}(u_0)\backslash\{u_1\}$, then $v_0$ must be a leaf vertex. By Lemma \ref{lem2.6}, we have $v_0\in I(\hat{T}),u_0\not\in I(\hat{T})$. Similarly, we also obtain $u_k\not\in I(\hat{T})$.

We firstly prove that $u_1\in I(\hat{T})$. Assume that $u_1\not\in I(\hat{T})$. Note that $d_{\hat{T}}(u_0)\geq3$ and $d_{\hat{T}}(u_k)\geq3$. Then the edge $u_0u_1$ is on an internal path of $\hat{T}$, which is contained in $P$. Denote by $T_1$ the tree built from $\hat{T}$ through subdividing the internal edge $u_0u_1$ and deleting a leaf vertex $v_0$ simultaneously. Note that $u_1\not\in I(\hat{T})$. It is straightforward to verify that $\gamma(T_1)=\gamma(\hat{T})$. By virtue of Lemmas \ref{lem2.1} and \ref{lem2.7}, one sees $\lambda_\alpha(T_1)<\lambda_\alpha(\hat{T})$, a contradiction to the minimality of $\lambda_\alpha(\hat{T})$.
Therefore, $u_1\in I(\hat{T})$ and $u_2\not\in I(\hat{T})$. 

Then we prove that $u_3\in I(\hat{T})$. Assume that $u_3\not\in I(\hat{T})$. Similarly, the edge $u_2u_3$ is also an internal edge of $\hat{T}$. Denote the tree $T_2$ constructed from $\hat{T}$ through subdividing the internal edge $u_2u_3$ and deleting a leaf vertex $v_0$ simultaneously. Clearly, $\gamma(T_2)=\gamma(\hat{T})$. By virtue of  Lemmas \ref{lem2.1} and \ref{lem2.7}, one gets $\lambda_\alpha(T_2)<\lambda_\alpha(\hat{T})$, which is also impossible.
Hence $u_3\in I(\hat{T})$. 

Using a similar method, we can show that $u_i\in I(\hat{T})$ when $i$ is odd, and $u_j\not\in I(\hat{T})$ when $j$ is even. Recall that $u_k\not\in I(\hat{T})$. Thus $k$ is even and $k\geq 2$, as desired. 
\end{proof}

\begin{lem}\label{Prop3.4}
Let $u_0,u_k$ be two end branch points of $\hat{T},$ and $P=u_0u_1
\cdots u_k(k\geq2$ $is~even)$ be a path that connects $u_0$ and $u_k$. Then for any vertex $u_i(i\in[1,k-1])$, it cannot attach a  pendant path of length at least two. Moreover, when $i$ is odd, $u_i$ cannot attach any pendant edge.
\end{lem}
\begin{proof}
Suppose, by way of contradiction, that there exists some $i\in[1,k-1]$, so that $u_i$ attaches a pendant path $P_1=u_iy_1\cdots y_s$ with $s\geq2$.
Note that $d_{\hat{T}}(u_0)\geq3$ and $d_{\hat{T}}(u_i)\geq3$.  Then the edge $u_{i-1}u_i$ is on an internal path of $\hat{T}$, which is contained in $P$.
Let $T_1$ be obtained from $\hat{T}$ by substituting the edge $u_{i-1}u_i$ by a path $u_{i-1}w_1w_2u_i$. Let $T_2=T_1-y_s-y_{s-1}$. By the same reasoning as that of the proof in Lemma \ref{Prop3.2}, we obtain 
$$|\hat{T}|=|T_1|-2=|T_2|,~\lambda_\alpha(\hat{T})\geq\lambda_\alpha(T_1)>\lambda_\alpha(T_2),~\gamma(\hat{T})=\gamma(T_1)-1=\gamma(T_2),$$
which contradicts the minimality of $\hat{T}$.

By virtue of Lemma \ref{Prop3.3}, one sees $u_i\in I(\hat{T})$ for odd $i$. This indicates that, when $i$ is odd, $u_i$ cannot attach any leaf vertex since $L(\hat{T})\subseteq I(\hat{T})$.
\end{proof}

\begin{lem}\label{Prop3.5}
Let $u_0,u_k$ be two end branch points of $\hat{T}$ and $P=u_0u_1\cdots u_k\, (k\geq2$ $is~even)$ be a path that connects $u_0$ and $u_k$. For odd $i\in[1,k-1]$, if $d_{\hat{T}}(u_i)\geq3$, then none of the vertices in $N_{\hat{T}}(u_i)$ attaches a leaf.
\end{lem}
\begin{proof}
Let $N_{\hat{T}}(u_i)=\{u_{i-1},u_{i+1},w_1,w_2,\dots,w_r\}$. Since $d_{\hat{T}}(u_i)\geq 3$, we have $r\geq1$. Suppose that there exists a vertex, say $w_1,$ in $N_{\hat{T}}(u_i)$ attaching a leaf, say $v$. Let $T_1=\hat{T}-u_iu_{i-1}+vu_{i-1}$. 
Obviously, $I(\hat{T})$ is still an independent set of $T_1$, which implies $\gamma(T_1)\geq \gamma(\hat{T})$. 
In what follows, we proceed by proving $\gamma(T_1)\leq\gamma(\hat{T})$. If $u_i\not\in I(T_1),$ then $I(T_1)$ is still an independent set of $\hat{T}$ and $\gamma(\hat{T})\geq\gamma(T_1)$. So in what follows, we consider $u_i\in I(T_1)$. In this case, we have $w_j\not\in I(T_1)$ for $1\leq j\leq r$ and $u_{i+1}\not\in I(T_1)$. Then we  differentiate the following two cases.

{\bf Case 1.} $v\in I(T_1)$. Clearly, $u_{i-1}\not\in I(T_1)$, which implies $I(T_1)$ is still an independent set of $\hat{T}$ and $\gamma(\hat{T})\geq\gamma(T_1)$. 

{\bf Case 2.} $v\not\in I(T_1)$. 
If $u_{i-1}\not\in I(T_1)$, then $I(T_1)$ is still an independent set of $\hat{T}$ and $\gamma(\hat{T})\geq\gamma(T_1)$. 
If $u_{i-1}\in I(T_1)$, let $S(\hat{T})=I(T_1)-u_{i-1}+v$, then $S(\hat{T})$ is an independent set of $\hat{T}$ satisfying  $|S(\hat{T})|=|I(T_1)|$, which implies $\gamma(\hat{T})\geq\gamma(T_1)$. 

Hence, we obtain $\gamma(T_1)=\gamma(\hat{T})$, $|T_1|=|\hat{T}|$. Then $T_1\in\mathcal{G}_{n,\gamma}$. 
Assume $\bf{x}$ is $A_\alpha(T_1)$'s Perron vector. By virtue of Lemma \ref{lem2.2}, 
if $x_{u_i}\geq x_{v}$, then $\lambda_\alpha(T_1)<\lambda_\alpha(T_1-vu_{i-1}+u_iu_{i-1})=\lambda_\alpha(\hat{T}).$
Similarly, if $x_{u_i}< x_{v}$, then $\lambda_\alpha(T_1)<\lambda_\alpha(T_1-u_iu_{i+1}-\{u_iw_s|s=2,\dots,r\}+vu_{i+1}
+\{vw_s|s=2,\dots,r\})=\lambda_\alpha(\hat{T}).$

Thus, one gets $\lambda_\alpha(T_1)<\lambda_\alpha(\hat{T})$, which will not happen.
\end{proof}

\begin{lem}\label{Prop3.7}
For any vertex $x\in V(\hat{T})\setminus L(\hat{T})$, there exist two end branch points $u_0,u_k$, such that $x$ is on the path $P=u_0u_1\cdots u_k\, (k\geq2$ $is~even)$, namely $x=u_i$ for some $0\leq i\leq k$. Furthermore, the labeling of $x=u_i$ is determined by the choice of path $P$, while the parity of its subscript $i$ remains independent of $P$.
\end{lem}
\begin{proof}
Suppose there is a vertex $x\in V(\hat{T})\setminus L(\hat{T})$ which is not on the path connecting any two end branch points of $\hat{T}$. Clearly, $d_{\hat{T}}(x)\geq 2$. Next, we continue by taking into account the following two possible scenarios.

{\bf Case 1.} $d_{\hat{T}}(x)\geq 3$. Then, $x$ is a branch point of $\hat{T}$. 

It suffices to prove that $x$ is not an end branch point. Suppose not, then there exists another end branch point $y$ (based on Lemma \ref{Prop3.1}), such that $x, y$ are connected in $\hat{T}$. Thus $x$ lies on the path between two end branch points in $\hat{T}$, which contradicts the assumption. Then $x$ is not an end branch point. 

By the definition of the end branch point, there exist other two branch point $u$ and $v$, so that $x$ lies on the path between  $u$ and $v$. It is straightforward to verify that $u$ and $v$ belong to the distinct components of $\hat{T}-x$, respectively. Denote the component containing $u$ by $T_{u}$, the component containing $v$ by $T_{v}$. 
Consider in $\hat{T}$, suppose that for any branch point in $V(T_u)$, it cannot be the end branch point in $\hat{T}$. Then for any branch point $u_i$ in $V(T_u)$, there exist other two branch points $u_j,u_k$, such that $u_i$ lies on the path joining $u_j$ and $u_k$ in $\hat{T}$. This leads to a contradiction since $\hat{T}$ is a tree and the number of the branch points in $\hat{T}$ is finite. Then there exists an branch point $u_1$ in $V(T_u)$, such that $u_1$ is an end branch point of $\hat{T}$. Similarly, there exists an branch point $v_1$ in $V(T_v)$, such that $v_1$ is an end branch point in $\hat{T}$. Thus $x$ lies on the path between two end branch points $u_1, v_1$ in $\hat{T}$, which contradicts the assumption.

{\bf Case 2.} $d_{\hat{T}}(x)=2$. Denote by $T_1, T_2$  the two components of $\hat{T}-x$, respectively.

{\bf Subcase 2.1.} There exist vertices $u_1\in V(T_1), v_1\in V(T_2)$, such that $d_{\hat{T}}(u_1)\geq 3$, $d_{\hat{T}}(v_1)\geq 3$. In view of the discussion of Case 1, similarly there must exist vertices $u_2\in V(T_1)$, $v_2\in V(T_2)$, such that $v_2$, $u_2$ are end branch points in $\hat{T}$. Then $x$ lies on the path between two end branch points $u_2, v_2$ in $\hat{T}$, which contradicts the assumption.

{\bf Subcase 2.2.} Let $w$ be in $V(\hat{T})\setminus\{x\}$. Then $d_{\hat{T}}(w)\leq2$. Clearly, $\hat{T}\cong P_n$, a contradiction since ${\lceil\frac{n}{2}\rceil}+2\leq\gamma(\hat{T})\leq n-2$.

{\bf Subcase 2.3.} Let $w$ be in $V(T_1)$. Then we have $d_{\hat{T}}(w)\leq2$, and there is a $y\in V(T_2)$ so that $d_{\hat{T}}(y)\geq 3$. Thus, $x$ attaches a path of length $t(t\geq1)$. Denote the branch point which is closest to $x$ by $v_1$. Then for any vertex $z$ on the path between $x$ and $v_1$, we have $d_{\hat{T}}(z)=2$. It is straightforward to check that $v_1$ attaches a path of length $s\, (s\geq 2)$. 

According to Lemma \ref{Prop3.2}, for any end branch point $u_0$ of $\hat{T}$, it attaches exactly $d_{\hat{T}}(u_0)-1$ leaves. It implies that $v_1$ can not be an end branch point of $\hat{T}$. By means of a discussion analogous to that of Case 1, there exist two end branch points $w_1, w_2$, such that $v_1$ is on the path between $w_1$ and $w_2$. By Lemma \ref{Prop3.4}, $v_1$ can not attach a pendant path of length at least 2, which leads to a contradiction.
Therefore, for any vertex $x\in V(\hat{T})\setminus L(\hat{T})$, there exist two end branch points $u_0,u_k$, such that $x$ is on the path $P=u_0u_1\cdots u_k\, (k\geq2$ $is~even)$, i.e., $x=u_i$ for some $0\leq i\leq k$.

Assume that there exist other two end branch points $v_0, v_s$, such that $x$ is also on the path $P^1=v_0v_1\cdots v_s\, (s\geq2$ is even) and $x=v_j$. Suppose to the contrary that $i$ and $j$ has different parity. For convenience, let $i$ be odd and $j$ be even. Then $d_{\hat{T}}(u_0,x)$ is odd, $d_{\hat{T}}(v_0,x)$ is even. Obviously, $x$ is on the path $P^2$ between two end branch points $u_0, v_0$ and the length of $P^2$ is odd, which contradicts the conclusion of Lemma \ref{Prop3.3}.
\end{proof}

According to Lemmas \ref{Prop3.3} and \ref{Prop3.7}, it is routine to determine the independent set $I(\hat{T})$ that contains all leaves of $\hat{T}$.  Let $V_1^*:=\{u_i\in V(\hat{T})\setminus L(\hat{T})|\, i$ is odd\}, and $V_2^*:=\{u_i\in V(\hat{T})\setminus L(\hat{T})|\, i$ is even\}. Then $I(\hat{T})$ is the disjoint union of $L(\hat{T})$ and $V_1^*$, and let $\bar{I}(\hat{T})=V(\hat{T})\setminus I(\hat{T})=V_2^*$.
\begin{lem}\label{Prop3.6}
$\hat{T}$ contains at least $2 \gamma-n+1$ leaf vertices.
\end{lem}
\begin{proof}
To prove Lemma \ref{Prop3.6}, it is equivalent to show that the subgraph $\hat{T}-L(\hat{T})$ has at most $2(n-\gamma)-1$ vertices.  We characterize the structure of $\hat{T}-L(\hat{T})$ in what follows.

Deleting the leaf vertices of $\hat{T}$ yields a non-trivial tree, and so it has at least two leaves. In view of Lemma \ref{Prop3.7},  for any vertex $x\in V(\hat{T})\setminus L(\hat{T})$, there exist two end branch points $u_0,u_k$, such that $x$ is on the path $P=u_0u_1\cdots u_k$, i.e., $x=u_i$ for some $0\leq i\leq k$. Obviously, we have $d_{\hat{T}-L(\hat{T})}(u_i)\geq 2$ if $i\in[1,k-1],$ and $d_{\hat{T}-L(\hat{T})}(u_i)= 1$ if $i\in\{0,k\}.$
Therefore, any leaf vertex in $\hat{T}-L(\hat{T})$ must be an end branch point of $\hat{T}$.

Consider any two distinct leaf vertices, say $u,v$, in $\hat{T}-L(\hat{T})$, the unique path connecting $u, v$ in $\hat{T}-L(\hat{T})$ is denoted by $P_{u,v}$. Obviously, $P_{u,v}$ is also the path that connects $u$ and $v$ in $\hat{T}$. In view of Lemma \ref{Prop3.3}, $u,v$ not in $I(\hat{T})$, and all the vertices on $P_{u,v}$ appear alternately in $I(\hat{T})$. Let $\bar{I}(\hat{T})=V(\hat{T})\setminus I(\hat{T})$. Then all the vertices on $P_{u,v}$ appear alternately in $\bar{I}(\hat{T})$ and $u,v\in \bar{I}(\hat{T})$. Moreover, $\bar{I}(\hat{T})$ contains all leaves of $\hat{T}-L(\hat{T})$.

Recall that $L(\hat{T})\subseteq I(\hat{T})$. Then $I(\hat{T})\setminus L(\hat{T})$ forms an independent set belonging to $\hat{T}-L(\hat{T})$. Additionally, for all $u,\,u'\in \bar{I}(\hat{T}),$ they cannot be adjacent in $\hat{T}$. Otherwise, $uu'$ must be on some path between two end branch points. It implies that $\bar{I}(\hat{T})$ is an independent set of 
$\hat{T}-L(\hat{T})$.
Therefore, $\hat{T}-L(\hat{T})$ has a bipartition $(I(\hat{T})\setminus L(\hat{T}), \bar{I}(\hat{T}))$.

Let $L_1$ be the set of leaves of $\hat{T}-\ L(\hat{T})$, $I_1$ be the maximum independent set of $\hat{T}- L(\hat{T})$ such that $L_1\subseteq I_1$. In what follows, we proceed by proving $I_1=\bar{I}(\hat{T})$.

Suppose, for the sake of contradiction, that $I_1\not=\bar{I}(\hat{T})$. Let $I_1$ be the disjoint union of $V_1, V_2$ and $L_1$, where $\emptyset \not =V_1\subseteq I(\hat{T})\setminus L(\hat{T})$, $V_2\subseteq \bar{I}(\hat{T})\setminus L_1$. In view of Lemmas \ref{Prop3.3} and \ref{Prop3.4}, for all $v\in I(\hat{T})\setminus L(\hat{T})$, one sees 
$$d_{\hat{T}-L(\hat{T})}(v)=d_{\hat{T}}(v)\geqslant 2, \ ~\mbox{and}~\ N_{\hat{T}}(v)=N_{\hat{T}-L(\hat{T})}(v)\subseteq \bar{I}(\hat{T}).$$ 
Let $V_3=\{N_{\hat{T}-L(\hat{T})}(u)| u\in V_1\}$. Then $V_3\subseteq \bar{I}(\hat{T})\setminus L_1$ and $|V_3|\geq |V_1|+1$. Moreover, for any vertex in $V_3$, it does not belong to $V_2$. Then $|V_2|\leq |\bar{I}(\hat{T})\setminus L_1|-(|V_1|+1)$. Thus, we obtain
$$|I_1|=|V_1|+|V_2|+|L_1|\leq |V_1|+|\bar{I}(\hat{T})\setminus L_1|-(|V_1|+1)+|L_1|\leq |\bar{I}(\hat{T})|-1,$$
a contradiction.

Therefore,  $\hat{T}- L(\hat{T})$ contains a maximum independent set $\bar{I}(\hat{T})$, which contains all leaves of $\hat{T}- L(\hat{T})$. Together with Lemma \ref{lem3.7}, $|\hat{T}- L(\hat{T})|\leqslant 2|\bar{I}(\hat{T})|-1=2(n-\gamma)-1$, as desired.
\end{proof}

\subsection{\normalsize Proof of Theorem~\ref{thm4.1}}
To accomplish the proof pertaining to Theorem~\ref{thm4.1}, we necessitate certain preliminaries. Denote by 
\[t =\left\lfloor\frac{2\gamma-n+1}{n-\gamma}\right\rfloor.\label{4.01}
\]

\begin{lem}\label{lem4.5}
Let $\hat{T}$ be in $\mathcal{G}_{n,\gamma}$ having minimal $A_\alpha$ spectral radius  
with $\gamma\in[\lceil\frac{n}{2}\rceil+2, n-2]$ and $\alpha\in[\frac{1}{2},1)$.
Then
$\lambda_\alpha(\hat{T})<\alpha( t +3)+2(1-\alpha).$
Particularly, if $|L(\hat{T})|\leq(n-\gamma) t +2$, then
$\lambda_\alpha(\hat{T})\leq\alpha( t +2)+2(1-\alpha).$
\end{lem}
\begin{proof}
Let $Y=\{v_0,v_2,\dots,v_{2(n-\gamma)-2}\}$ be the rooted set of the path $P_{2(n-\gamma)-1}= v_0v_1\cdots
v_{2(n-\gamma)-2}$, and let $T=P_{2(n-\gamma)-1}\circ_Y(t +1)$. By virtue of Lemma \ref{lem4.2.1}, we have
$\lambda_\alpha(T)\leq\alpha( t +3)+2(1-\alpha).$
Obviously, $|L(T)|=( t +1)|Y|=(\lfloor\frac{2\gamma-n+1}{n-\gamma}\rfloor+1)(n-\gamma)>2\gamma-n+1\geq5$. Let \(T_1\) be the graph derived from $T$ by removing some leaves, with the condition that $|L(T_1)|=2\gamma-n+1$ and $v_0$, $v_{2(n-\gamma)-2}$ each have at least one leaf vertex attached in $T_1$. Recall that the maximum independent set $I(T_1)$ includes all the leaf vertices of $T_1$. Then $v_0\not\in I(T_1), v_{2(n-\gamma)-2}\not\in I(T_1)$ and we obtain
$$|I(T_1)|=\gamma(T_1)\leq\Big\lceil\frac{2(n-\gamma)-1-2}{2}\Big\rceil+(2\gamma-n+1)=\gamma.$$ 
Observe that  $L(T_1)\cup (V(P_{2(n-\gamma)-1})\setminus Y)$ is a disjoint union, which forms an independent set of $T_1$, and $|L(T_1)|+(|P_{2(n-\gamma)-1}|-|Y|)=\gamma$.
Therefore, we have $\gamma(T_1)=\gamma$, $|T_1|=|L(T_1)|+|P_{2(n-\gamma)-1}|=n$. That is, $T_1\in \mathcal{G}_{n,\gamma}$. Recall that $T_1$ is properly contained in  $T$. In view of Lemma \ref{lem2.1}, one obtains $\lambda_\alpha(\hat{T})\leq\lambda_\alpha(T_1)<\lambda_\alpha(T)$.

Next one may consider $|L(\hat{T})|\leq(n-\gamma) t +2$.
Let $T'=P_{2(n-\gamma)-1}\circ_Y( t +1, t , t ,\dots, t , t , t +1)$.
By Lemma \ref{lem4.2.1}, one has 
$\lambda_\alpha(T')\leq\alpha( t +2)+2(1-\alpha).$
According to Lemma \ref{Prop3.6}, $|L(T')|= 2+t |Y|=(n-\gamma) t +2\geq |L(\hat{T})|\geq2\gamma-n+1\geq5$.
Deleting some leaves from $T'$, we can obtain a tree $T_2$, such that $L(T_2)=2\gamma-n+1$ and $v_0$, $v_{2(n-\gamma)-2}$ attach at least one leaf vertex in $T_2$, respectively. By a similar discussion, we obtain $T_2\in \mathcal{G}_{n,\gamma}$. Notice that $T_2$ is a subgraph of $T'$. In light of Lemma \ref{lem2.1}, one has $\lambda_\alpha(\hat{T})\leq\lambda_\alpha(T_2)\leq\lambda_\alpha(T')$, as required.
\end{proof}

\begin{lem}\label{lem4.4}
Let $\hat{T}$ be in $\mathcal{G}_{n,\gamma}$ having the minimal $A_\alpha$ spectral radius, where $\gamma\in[\lceil\frac{n}{2}\rceil+2, n-2]$ and $\alpha\in[0,\frac{1}{2})$.
Then 
$\lambda_\alpha(\hat{T})<\alpha( t +5)+(1-2\alpha)\sqrt{ t +5}$, where $t$ is defined in \eqref{4.01}
\end{lem}
\begin{proof}
Let $Y=\{v_0,v_2,\dots,v_{2(n-\gamma)-2}\}$ be the rooted set of the path $P_{2(n-\gamma)-1}= v_0v_1\cdots
v_{2(n-\gamma)-2}$, and let $T=P_{2(n-\gamma)-1}\circ_Y( t +1)$. 
In light of Lemma \ref{lem4.1}, one sees 
\begin{equation}\label{4.1}
\lambda_\alpha(T)\leq2\alpha\lambda_{\frac{1}{2}}(T)+(1-2\alpha)\lambda_{0}(T).
\end{equation}
Next we may consider the values of $\lambda_{\frac{1}{2}}(T)$ and $\lambda_{0}(T)$, respectively.

According to the proof of Lemma \ref{lem4.5}, one gets $\lambda_{\frac{1}{2}}(T)\leq\frac{1}{2}( t +5)$. By Lemma \ref{lem4.0}, we have 
$$\lambda_{0}(T)=\sqrt{{\lambda_{0}}^2(P_{2(n-\gamma)-1})+ t +1}.$$
Note that $\lambda_{0}(P_{2(n-\gamma)-1})=2\cos(\frac{\pi}{2(n-\gamma)})<2$ (based on Lemma \ref{lem4.2.2}). Thus, we obtain $\lambda_{0}(T)<\sqrt{ t +5}$. Together with \eqref{4.1}, we have
$$\lambda_\alpha(T)<\alpha( t +5)+(1-2\alpha)\sqrt{ t +5}.$$

By the same discussion as that of Lemma \ref{lem4.5}, one may obtain a graph $T_1$ from $T$ by removing some leaf vertices, such that  $T_1\in\mathcal{G}_{n,\gamma}$. By virtue of Lemma \ref{lem2.1}, one gets $\lambda_\alpha(\hat{T})\leq\lambda_\alpha(T_1)<\lambda_\alpha(T)$, as required.
\end{proof}

\begin{lem}\label{lem4.7}
Assume $c\geq4$ is a given integer, $\gamma=(n-c)\in[\lceil\frac{n}{2}\rceil+2,n-4]$ and $\hat{T}$ is an $A_\alpha$-minimizer graph in $\mathcal{G}_{n,\gamma}$ with $0\leq\alpha < 1$. Let $u_0,u_k$ be two end branch points of $\hat{T}$ and $P=u_0u_1\cdots u_k\, (k\geq2$ $is~even)$ be a path that connects $u_0$ and $u_k$. Then \(d_{T^*}(u_i) = 2\) for each odd \(i \in [1, k-1]\) if any of the following conditions is met.
\begin{enumerate}[{\rm(i)}]
    \item $\alpha=0, n\geq \lceil\frac{4c^2-6c-3}{3}\rceil$.
    \item $\alpha\in(0,\frac{1}{2}), n\geq \max\{{c(\frac{(1-2\alpha)(c-3)}{3\alpha}+\sqrt{c+2})^2-3c-1}, 2c+4\}$.
    \item $\alpha\in[\frac{1}{2},1), n\geq \max\{c^2-c-1, 2c+4\}$. 
\end{enumerate}
\end{lem}
\begin{proof}
Notice that $\gamma=(n-c)\in[\lceil\frac{n}{2}\rceil+2,n-4]$. Hence, we have $n\geq 2c+4$.
Moreover, the $A_\alpha$-minimizer graph $\hat{T}$ satisfies Lemmas \ref{Prop3.1}-\ref{Prop3.6}. 

Suppose that there exists some odd $i\in [1,k-1]$ such that $d_{\hat{T}}(u_i)\geq3$ for $\alpha\in[0,1)$. According to Lemma \ref{Prop3.5}, for any vertex in $N_{\hat{T}}(u_i)$, it cannot attach a leaf in $\hat{T}$. Let $\bar{I}(\hat{T})=V(\hat{T})\setminus I(\hat{T})$, then any leaf vertex must be a neighbour of some vertex of  $\bar{I}(\hat{T})$ since $I(\hat{T})$ contains all leaves of $\hat{T}.$ In light of Lemma \ref{Prop3.3}, one sees $u_i\in I(\hat{T})$ and $N_{\hat{T}}(u_i)\subseteq \bar{I}(\hat{T})$. Hence, $\vert \cup_{x \in L(\hat{T})} N_{\hat{T}}(x)\vert \leqslant |\bar{I}(\hat{T})|-|N_{\hat{T}}(u_i)|=n-\gamma-d_{\hat{T}}(u_i)$. By virtue of  Lemma \ref{Prop3.6}, one sees $\hat{T}$ has at least $2\gamma-n+1$ leaf vertices. Thus, one can find a vertex $u$ which is adjacent to at least $\lceil\frac{2\gamma-n+1}{n-\gamma-d_{\hat{T}}(u_i)}\rceil$ leaf vertices. Since $d_{\hat{T}}(u_i)\geq3$, we have $\frac{2\gamma-n+1}{n-\gamma-d_{\hat{T}}(u_i)}\geq\frac{2\gamma-n+1}{n-\gamma-3}$. Moreover, we obtain
$$d_{\hat{T}}(u)\geqslant\lceil\frac{2\gamma-n+1}{n-\gamma-d_{\hat{T}}(u_i)}\rceil+1\geqslant\frac{2\gamma-n+1}{n-\gamma-3}+1=\frac{\gamma-2}{n-\gamma-3}.$$
Note that $\hat{T}$ contains $K_{1,d_{\hat{T}}(u)}$ as a proper subgraph. Along with Lemma \ref{lem2.1}, we get  $\lambda_\alpha(\hat{T})>\lambda_\alpha(K_{1,d_{\hat{T}}(u)}).$
Below, let $ \bar{\ell} =\frac{2\gamma-n+1}{n-\gamma}$. 

When $\alpha=0$, since $n\geq \lceil\frac{4c^2-6c-3}{3}\rceil$, it is routine to check that $n>2c+4$. By Lemma \ref{lem4.2}, we have
$$\lambda_0(\hat{T})>\lambda_0(K_{1,d_{\hat{T}}(u)})=\sqrt {d_{\hat{T}}(u)}\geq \sqrt{\frac{\gamma-2}{n-\gamma-3}}
\geq \sqrt{ \bar{\ell} +5},$$
which contradicts Lemma \ref{lem4.4}.

When $\alpha\in(0,\frac{1}{2})$, since $n\geq\max\{{c(\frac{(1-2\alpha)(c-3)}{3\alpha}+\sqrt{c+2})^2-3c-1}, 2c+4\}$, in light of Lemma \ref{lem4.2}, one gets
$$\lambda_\alpha(\hat{T})>\lambda_\alpha(K_{1,d_{\hat{T}}(u)})\geq \alpha(d_{\hat{T}}(u)+1)\geq \alpha(\frac{\gamma-2}{n-\gamma-3}+1)
\geq \alpha(\bar{\ell}+5)+(1-2\alpha)\sqrt{ \bar{\ell} +5},$$
which contradicts Lemma \ref{lem4.4}.

When $\alpha\in[\frac{1}{2},1)$, since $n\geq \max\{c^2-c-1, 2c+4\}$, in light of Lemma \ref{lem4.2}, one has 
$$\lambda_\alpha(\hat{T})>\lambda_\alpha(K_{1,d_{\hat{T}}(u)})\geq \alpha d_{\hat{T}}(u)+\frac{(1-\alpha)^2}{\alpha}\geq \alpha\frac{\gamma-2}{n-\gamma-3}+\frac{(1-\alpha)^2}{\alpha}
\geq \alpha( \bar{\ell} +3)+2(1-\alpha),$$
which contradicts Lemma \ref{lem4.5}.
\end{proof}

Now, it is ready to show Theorem~\ref{thm4.1} armed with these results.
\begin{proof}[\bf Proof of Theorem \ref{thm4.1}]
With the help of the proof of Lemma \ref{Prop3.6}, one sees $\hat{T}- L(\hat{T})$ is a tree, which has a bipartition $(I(\hat{T})\setminus L(\hat{T}), \bar{I}(\hat{T}))$ with $\bar{I}(\hat{T})=V(\hat{T})\setminus I(\hat{T})$. That is, generality is not lost if we assume that, for any edge $uv\in E(\hat{T}-L(\hat{T}))$, $u$ is in $I(\hat{T})\setminus L(\hat{T})$, and $v$ is in $\bar{I}(\hat{T})$.
In view of Lemmas \ref{Prop3.3} and \ref{lem4.7}, every vertex of $\hat{T}-L(\hat{T})$ in $I(\hat{T})\setminus L(\hat{T})$ is of degree 2. By Lemma \ref{lem4.3}, $\hat{T}-L(\hat{T})$ is constructed from some tree $T$ by subdividing every edge of $T$, where $|T|$=$|\bar{I}(\hat{T})|=n-\gamma$. Thus $|\hat{T}-L(\hat{T})|=2(n-\gamma)-1$. Clearly, $|L(\hat{T})|=n-(2(n-\gamma)-1)=2\gamma-n+1$. Recall that $L(\hat{T})\subseteq I(\hat{T})$, then for any leaf of $\hat{T}$, it can only be attached to some vertex in $\bar{I}(\hat{T})=V(T)$.
This completes the proof.  
\end{proof}
\section{\normalsize Concluding remarks}\setcounter{equation}{0}
Throughout this section, we aim to characterize the graph in $\mathcal{G}_{n,n-4}$ having minimal $A_\alpha$ spectral radius, where $\alpha\in[\frac{1}{2},1)$ and $n\geq 12$. 
It can be seen as a continuince of the work of Liu and Wang \cite{13}, in which they characterized the $A_\alpha$-minimizer graphs in $\mathcal{G}_{n,n-4}$ for $\alpha\in[\frac{1}{2},1)$. Together with these results, we identify all the possible $A_\alpha$-minimizer graphs, which extends the results obtained in \cite{13}.

Recall that 
$t =\left\lfloor\frac{2\gamma-n+1}{n-\gamma}\right\rfloor$. Then let $\ell'=(2\gamma-n+1)-(n-\gamma) t.$
Clearly, $0\leq\ell'\leq n-\gamma-1$.

\begin{lem}\label{lem4.6}
Let $\hat{T}$ be in  $\mathcal{G}_{n,\gamma}$ having minimal \(A_\alpha\) spectral radius. Let $\ell(u)$ represent the count of leaf vertices attached to a vertex $u\in V(T)$.
\begin{enumerate}[{\rm(a)}]
    \item If~$\alpha\in (0,\frac{1}{4})$, then
    $0\leq\ell(u)\leq\lfloor\frac{8}{5}( t +5)\rfloor-d_{T}(u).$
   \item If~$\alpha\in [\frac{1}{4},\frac{1}{2})$, then 
   $0\leq\ell(u)\leq\lfloor(4-2\sqrt{2})( t +5)\rfloor-1-d_{T}(u).$
    \item If~$\alpha\in [\frac{1}{2},1)$, then
    \begin{enumerate}[{\rm(i)}]
      \item  
      $ t +\ell'-d_{T}(u)\leq\ell(u)\leq t +2-d_{T}(u)$ for $0\leq\ell'\leq2$.
      \item 
      $ t +\ell'-(n-\gamma-1)-d_{T}(u)\leq\ell(u)\leq t +3-d_{T}(u)$ for $3\leq\ell'\leq n-\gamma-1$.
    \end{enumerate}
\end{enumerate}
\end{lem}
\begin{proof}
(a)\ Notice that $K_{1,d_{\hat{T}}(u)}$ is properly contained in $\hat{T}$ for all vertices $u$ in $V(T)$. In view of Lemmas \ref{lem2.1}, \ref{lem4.2} and \ref{lem4.4}, one gets 
\begin{align}\label{4.2}
\lambda_\alpha(K_{1,d_{\hat{T}}(u)})&<\lambda_\alpha(\hat{T}),\\
\lambda_\alpha(\hat{T})&<\alpha( t +5)+(1-2\alpha)\sqrt{ t +5},\notag\\
\lambda_\alpha(K_{1,d_{\hat{T}}(u)})&=\frac{1}{2}(\alpha(d_{\hat{T}}(u)+1)+\sqrt{\alpha^2(d_{\hat{T}}(u)+1)^2+4d_{\hat{T}}(u)(1-2\alpha)}).\notag
\end{align}
Assume that $d_{\hat{T}}(u)\geq\lceil\frac{8}{5}( t +5)\rceil$, then by direct calculations we obtain $$\lambda_\alpha(K_{1,d_{\hat{T}}(u)})>\alpha( t +5)+(1-2\alpha)\sqrt{ t +5}>\lambda_\alpha(\hat{T}),$$ which contradicts \eqref{4.2}. Thus, we have $d_{\hat{T}}(u)\leq\lfloor\frac{8}{5}( t +5)\rfloor$. Moreover,
$$0\leq\ell(u)=d_{\hat{T}}(u)-d_{\hat{T}-L(\hat{T})}(u)=d_{\hat{T}}(u)-d_{T}(u)\leq\lfloor\frac{8}{5}( t +5)\rfloor-d_{T}(u).$$

(b)\ By Lemma \ref{lem4.2}, we have
\begin{align*}
   \lambda_\alpha(K_{1,d_{\hat{T}}(u)})=&\frac{1}{2}(\alpha(d_{\hat{T}}(u)+1)+\sqrt{\alpha^2(d_{\hat{T}}(u)+1)^2+4d_{\hat{T}}(u)(1-2\alpha)})\\
   \geq&\frac{2+\sqrt2}{4}\alpha(d_{\hat{T}}(u)+1)+\frac{\sqrt{2-4\alpha}}{2}\sqrt{d_{\hat{T}}(u)}.
\end{align*}
Assume that $d_{\hat{T}}(u)\geq \lceil(4-2\sqrt{2})( t +5)\rceil-1$, then it is straightforward to check that $\frac{2+\sqrt2}{4}(d_{\hat{T}}(u)+1)> t +5$. By direct calculations, we obtain $\frac{\sqrt{2-4\alpha}}{2}\geq 1-2\alpha$ for $\alpha\in[\frac{1}{4}, \frac{1}{2})$. Therefore, when $d_{\hat{T}}(u)\geq \lceil(4-2\sqrt{2})( t +5)\rceil-1$, we obtain
\begin{align*}
   \lambda_\alpha(K_{1,d_{\hat{T}}(u)})
   \geq&\frac{2+\sqrt2}{4}\alpha(d_{\hat{T}}(u)+1)+\frac{\sqrt{2-4\alpha}}{2}\sqrt{d_{\hat{T}}(u)}\\
   >&\alpha( t +5)+(1-2\alpha)\sqrt{ t +5}\\
   >&\lambda_\alpha(\hat{T}),
\end{align*}
which contradicts \eqref{4.2}.  Thus, we have $d_{\hat{T}}(u)\leq\lfloor(4-2\sqrt{2})( t +5)\rfloor-1$. Moreover,
$$0\leq\ell(u)=d_{\hat{T}}(u)-d_{\hat{T}-L(\hat{T})}(u)=d_{\hat{T}}(u)-d_{T}(u)\leq\lfloor(4-2\sqrt{2})( t +5)\rfloor-1-d_{T}(u).$$

(c)\ If $0\leq\ell'\leq2$, then one sees $L(\hat{T})=(n-\gamma) t +\ell'\leq(n-\gamma) t +2$. In light of Lemma \ref{lem4.5}, one gets $\lambda_\alpha(\hat{T})\leq\alpha( t +2)+2(1-\alpha).$
Notice that $K_{1,d_{\hat{T}}(u)}$ is properly contained in $\hat{T}$ for all $u\in V(T)$. Together with Lemmas \ref{lem2.1} and \ref{lem4.2} gives us 
$$\alpha d_{\hat{T}}(u)+\frac{(1-\alpha)^2}{\alpha}\leq\lambda_\alpha(K_{1,d_{\hat{T}}(u)})<\lambda_\alpha(\hat{T})
\leq\alpha( t +2)+2(1-\alpha).$$
Thus, $d_{\hat{T}}(u)< t +2+\frac{2(1-\alpha)}{\alpha}-\frac{(1-\alpha)^2}{\alpha^2}\leq t +3$, that is, $d_{\hat{T}}(u)\leq t +2$. Moreover,
$$\ell(u)=d_{\hat{T}}(u)-d_{\hat{T}-L(\hat{T})}(u)=d_{\hat{T}}(u)-d_{T}(u)\leq t +2-d_{T}(u).$$
On the other hand, one has 
\begin{align*}
   \ell(u)=&|L(\hat{T})|-\sum_{\stackrel{v\in V(T)}{v\neq u}}\ell(v)=(n-\gamma) t +\ell'-\sum_{\stackrel{v\in V(T)}{v\neq u}}\ell(v)\\
   \geq&(n-\gamma) t +\ell'-\sum_{\stackrel{v\in V(T)}{v\neq u}}\Big( t +2-d_{T}(v)\Big)\\
   =&\ell'+(n-\gamma) t -(n-\gamma-1) t -2(n-\gamma-1)
   +\sum_{v\in V(T)}d_{T}(v)-d_{T}(u)\\
   =& t +\ell'-2(n-\gamma-1)+2|E(T)|-d_{T}(u)\\
   =& t +\ell'-d_{T}(u).
\end{align*}
Therefore, for $0\leq\ell'\leq2$, we obtain $ t +\ell'-d_{T}(u)\leq\ell(u)\leq t +2-d_{T}(u).$

If $3\leq\ell'\leq n-\gamma-1$, then by the same argument as that of the case $0\leq\ell'\leq2$, we can easily get 
$$ t +\ell'-(n-\gamma-1)-d_{T}(u)\leq\ell(u)\leq t +3-d_{T}(u),$$
which is based on Lemmas \ref{lem2.1}, \ref{lem4.2} and \ref{lem4.5}.

This completes the proof.
\end{proof}

\begin{remark}
\textup{According to Lemma \ref{lem4.6}, if $\alpha= \frac{1}{2}$, then} 
\textup{the range of $\ell(u)$ above concides with \cite[Lemma 4.2]{12}, which characterizes the $A_\frac{1}{2}$-minimizer graph.}
\end{remark}

\begin{thm}\label{thm 1.4}
Assume that $\hat{T}$ is the graph within $\mathcal{G}_{n,n-4}$ having minimal $A_\alpha$ spectral radius  with $n\geq 12$ and $\alpha\in[\frac{1}{2},1)$. Then $\hat{T}\cong S(K_{1,3})\circ_{V(K_{1,3})}(\ell(w_1),\ell(w_2),\ell(w_3),\ell(w_4))$, or $S(P_4)\circ_{V(P_4)}(\ell(u_1),\ell(u_2),\ell(u_3),\ell(u_4)),$
where $\sum_{i=1}^{4}\ell(w_i)=\sum_{i=1}^{4}\ell(u_i)=2\gamma-n+1=n-7$~(see Fig. $\ref{Fig.4}$). Moreover, all the possible $A_\alpha$-minimizer graphs are listed in the Table $\ref{Tab 1}$.  
\end{thm}
\begin{proof}
The first part of our result is a direct consequence of Theorem \ref{thm4.1}. Clearly, $ t =\lfloor\frac{n-7}{4}\rfloor\geq 1$ and $\ell'=0,1,2,3$. By Lemma \ref{lem4.6}, we can determine value ranges of $\ell(w_i)$ and $\ell(u_i)\, (i\in{1,2,3,4})$, respectively. We continue by looking at the following two cases. 

{\bf Case 1.} $\hat{T}\cong S(K_{1,3})\circ_{V(K_{1,3})}(\ell(w_1),\ell(w_2),\ell(w_3),\ell(w_4))$.
Note that $d_{K_{1,3}}(w_1)=d_{K_{1,3}}(w_2)=d_{K_{1,3}}(w_3)=1$ and $d_{K_{1,3}}(w_4)=3$. Together with Lemma \ref{lem4.6}, we have the following.

$\bullet$ $n-7\equiv 0 \pmod{4}$. Then
$\ell(w_1),\ell(w_2),\ell(w_3)\in[ t -1, t +1]~\mbox {and}~\ell(w_4)\in[ t -3, t -1].$

$\bullet$ $n-7\equiv 1 \pmod{4}$. Then
$\ell(w_1),\ell(w_2),\ell(w_3)\in[ t , t +1]~\mbox {and}~\ell(w_4)\in[ t -2, t -1].$

$\bullet$  $n-7\equiv 2 \pmod{4}$. Then
$\ell(w_1),\ell(w_2),\ell(w_3)= t +1~\mbox {and}~\ell(w_4)= t -1.$

$\bullet$  $n-7\equiv 3 \pmod{4}$. Then
$\ell(w_1),\ell(w_2),\ell(w_3)\in[ t -1, t +2]~\mbox {and}~\ell(w_4)\in[ t -3, t ].$

{\bf Case 2.} $\hat{T}\cong S(P_4)\circ_{V(P_4)}(\ell(u_1),\ell(u_2),\ell(u_3),\ell(u_4))$.
Note that $d_{P_4}(u_1)=d_{P_4}(u_4)=1$ and $d_{P_4}(u_2)=d_{P_4}(u_3)=2$. Together with Lemma \ref{lem4.6}, we have the following.

$\bullet$ $n-7\equiv 0 \pmod{4}$. Then 
$\ell(u_1),\ell(u_4)\in[ t -1, t +1]~\mbox {and}~\ell(u_2),\ell(u_3)\in[ t -2, t ].$

$\bullet$ $n-7\equiv 1 \pmod{4}$. Then 
$\ell(u_1),\ell(u_4)\in[ t , t +1]~\mbox {and}~\ell(u_2),\ell(u_3)\in[ t -1, t ].$

$\bullet$  $n-7\equiv 2 \pmod{4}$. Then 
$\ell(u_1),\ell(u_4)= t +1~\mbox {and}~\ell(u_2),\ell(u_3)= t .$

$\bullet$  $n-7\equiv 3 \pmod{4}$. Then 
$\ell(u_1),\ell(u_4)\in[ t -1, t +2]~\mbox {and}~\ell(u_2),\ell(u_3)\in[ t -2, t +1].$

By Cases 1 and 2, all the possible $A_\alpha$-minimizer graphs are given in Table \ref{Tab 1}. 
\end{proof}
\begin{table}[h!]
\caption{All the possible $A_\alpha$-minimizer graphs in $\mathcal{G}_{n,n-4}$ obtained in Theorem \ref{thm 1.4}.}\label{Tab 1}
  \centering
  \renewcommand{\arraystretch}{1.5}
  \begin{tabular}{p{1cm}<{\centering}|p{5cm}<{\centering}| p{5cm}<{\centering}}
\hline
$\ell'$&~$(\ell(w_1),\ell(w_2),\ell(w_3),\ell(w_4))$&$(\ell(u_1),\ell(u_2),\ell(u_3),\ell(u_4))$\\
\hline
0&\makecell*[l]{$( t -1,~ t +1,~ t +1,~ t -1)$\\
{$( t ,~ t +1,~ t +1,~ t -2)$}\\
{$( t ,~ t ,~ t +1,~ t -1)$}\\
{$( t +1,~ t +1,~ t +1,~ t -3)$}}&
\makecell*[l]{$( t -1,~ t ,~ t ,~ t +1)$\\
{$( t ,~ t -1,~ t ,~ t +1)$}\\
{$( t ,~ t ,~ t -1,~ t +1)$}\\
$( t ,~ t ,~ t ,~ t )$\\
$( t +1,~ t -2,~ t ,~ t +1)$\\
{$( t +1,~ t -1,~ t -1,~ t +1)$}}\\
\hline
1&\makecell*[l]{{$( t ,~ t +1,~ t +1,~ t -1)$}\\
{$( t +1,~ t +1,~ t +1,~ t -2)$}}&
\makecell*[l]{{$( t ,~ t ,~ t ,~ t +1)$}\\
{$( t +1,~ t -1,~ t ,~ t +1)$}}\\
\hline
2&\makecell[l]{{$( t +1,~ t +1,~ t +1,~ t -1)$}}&
\makecell[l]{{$( t +1,~ t ,~ t ,~ t +1)$}}\\
\hline
3&\makecell*[l]{$( t -1,~ t +2,~ t +2,~ t )$\\
$( t ,~ t +2,~ t +2,~ t -1)$\\
$( t ,~ t +1,~ t +2,~ t )$\\
{$( t +1,~ t +2,~ t +2,~ t -2)$}\\
{$( t +1,~ t +1,~ t +2,~ t -1)$}\\
$( t +1,~ t +1,~ t +1,~ t )$\\
{$( t +2,~ t +2,~ t +2,~ t -3)$}}&
\makecell*[l]{$( t -1,~ t +1,~ t +1,~ t +2)$\\
$( t ,~ t ,~ t +1,~ t +2)$\\
$( t ,~ t +1,~ t ,~ t +2)$\\
$( t ,~ t +1,~ t +1,~ t +1)$\\
{$( t +1,~ t +1,~ t -1,~ t +2)$}\\
{$( t +1,~ t -1,~ t +1,~ t +2)$}\\
{$( t +1,~ t ,~ t ,~ t +2)$}\\
{$( t +1,~ t ,~ t +1,~ t +1)$}\\
$( t +2,~ t -2,~ t +1,~ t +2)$\\
{$( t +2,~ t -1,~ t ,~ t +2)$}}\\
\hline
\end{tabular}
\end{table}
In what follows, we will use another methods to exclude some graphs in Table \ref{Tab 1} in order to make the results more accurate. Initially, we take into account the case  $ t \geq3$, which implies that $n\geq 19$. For convenience, let $T_1(\ell(w_1),\ell(w_2),\ell(w_3),\ell(w_4))=S(K_{1,3})\circ_{V(K_{1,3})}(\ell(w_1),\ell(w_2),\ell(w_3),\ell(w_4))$, and $T_2(\ell(u_1),\ell(u_2),\ell(u_3),\ell(u_4))=S(P_4)\circ_{V(P_4)}(\ell(u_1),\ell(u_2),\ell(u_3),\ell(u_4))$. By Lemma \ref{lem5.1}, we may exclude some graphs in Table \ref{Tab 1}.

If $\ell'=0$, then one has 
\begin{align*}
\lambda_\alpha(T_1( t , t , t +1, t -1))
&<\lambda_\alpha(T_1( t -1, t +1, t +1, t -1)),\\
\lambda_\alpha(T_2( t , t , t , t ))
&<\lambda_\alpha(T_2( t -1, t , t , t +1)),\\
\lambda_\alpha(T_2( t +1, t -1, t -1, t +1))
&<\lambda_\alpha(T_2( t +1, t -2, t , t +1)).
\end{align*}

If $\ell'=3$, then we have
\begin{align*}
\lambda_\alpha(T_1( t +1, t +1, t +1, t ))
&<\lambda_\alpha(T_1( t , t +1, t +2, t ))
<\lambda_\alpha(T_1( t -1, t +2, t +2, t )),\\
\lambda_\alpha(T_1( t +1, t +1, t +2, t -1))
&<\lambda_\alpha(T_1( t , t +2, t +2, t -1)),\\
\lambda_\alpha(T_2( t , t +1, t +1, t +1))
&<\lambda_\alpha(T_2( t -1, t +1, t +1, t +2)),\\
\lambda_\alpha(T_2( t +2, t -1, t , t +2))
&<\lambda_\alpha(T_2( t +2, t -2, t +1, t +2)).
\end{align*}

Based on Lemmas \ref{lem2.2} and \ref{lem5.2}, we can further exclude some graphs in Table \ref{Tab 1} for $\ell'=3$:
\begin{align*}
\lambda_\alpha(T_2( t +1, t -1, t +1, t +2))
&<\lambda_\alpha(T_2( t , t , t +1, t +2)),\\
\lambda_\alpha(T_2( t +1, t , t , t +2))
&<\lambda_\alpha(T_2( t , t +1, t , t +2)),\\
\lambda_\alpha(T_2( t +1, t , t +1, t +1))
&<\lambda_\alpha(T_2( t , t +1, t +1, t +1)).
\end{align*}

\begin{lem}\label{lem5.6}
Let $T_1=T_2( t , t -1, t , t +1)$, $T_2=T_2( t , t , t -1, t +1)$. Then $\lambda_\alpha(T_1)>\lambda_\alpha(T_2)$ for $\alpha\in[\frac{1}{2},1)$.
\end{lem}
\begin{proof}
Let $W_1=N_{T_1}(u_1)\setminus \{w_1\}$, $W_2=N_{T_1}(u_2)\setminus \{w_1, w_2\}$, $W_3=N_{T_1}(u_3)\setminus \{w_2, w_3\}$, $W_4=N_{T_1}(u_4)\setminus\{w_3\}$. It is routine to check that $\pi_1=\{W_1\}\cup\{u_1\}\cup\{w_1\}\cup\{W_2\}\cup\{u_2\}\cup\{w_2\}\cup\{W_3\}\cup\{u_3\}\cup\{w_3\}\cup\{W_4\}\cup\{u_4\}$
is an equitable partition of $V(T_1)$, whose quotient matrix is written as $A_\alpha(T_1)_{\pi_1}$. 

Let $W_1'=N_{T_2}(u_1)\setminus \{w_1\}$, $W_2'=N_{T_2}(u_2)\setminus \{w_1, w_2\}$, $W_3'=N_{T_2}(u_3)\setminus \{w_2, w_3\}$, $W_4'=N_{T_2}(u_4)\setminus\{w_3\}$. Clearly, $\pi_2=\{W_1'\}\cup\{u_1\}\cup\{w_1\}\cup\{W_2'\}\cup\{u_2\}\cup\{w_2\}\cup\{W_3'\}\cup\{u_3\}\cup\{w_3\}\cup\{W_4'\}\cup\{u_4\}$
is an equitable partition of $V(T_2)$, whose quotient matrix is written as  $A_\alpha(T_2)_{\pi_2}$. 

Let $f_1(\alpha, t ,x)=\det(xI-A_\alpha(T_1)_{\pi_1})$, $f_2(\alpha, t ,x)=\det(xI-A_\alpha(T_2)_{\pi_2})$. Assume that $x^1_{\text{max}}$ and $x^2_{\text{max}}$ are the largest zeros corresponding to $f_1(\alpha, t ,x)$ and $f_2(\alpha, t ,x)$, respectively.  Notice that $A_\alpha(T_1)$ and $A_\alpha(T_2)$ are all real nonnegative matrices. By virtue of Lemma \ref{lem5.3}, we have $x_{\max}^{1}=\lambda_\alpha(T_1)$ and $x_{\max}^{2}=\lambda_\alpha(T_2)$. Note that $\lambda_\alpha(T_1)$ and $\lambda_\alpha(T_2)$ must be real. Then $x_{\max}^{1}$ and $x_{\max}^{2}$ are the largest real zeros of $f_1(\alpha, t ,x)$ and $f_2(\alpha, t ,x)$, respectively.
Together with Lemma \ref{lem4.2}, we have
\begin{equation}\label{5.3}
  \alpha( t +2)+\frac{(1-\alpha)^2}{\alpha}\leq\lambda_\alpha(T_2)=x_{\max}^{2}.
\end{equation}

Utilizing calculations by MATLAB R2023\cite{15}, we obtain
\begin{equation*}
  f_1(\alpha, t ,x)-f_2(\alpha, t ,x)=(\alpha-1)^4(2\alpha-x)[\alpha x^2-(\alpha^2+2\alpha-1)x+2\alpha^2-\alpha]^2.
\end{equation*}

It is routine to check that $x_1:=\alpha$, $x_2:=2-\frac{1}{\alpha}$ are the roots of $\alpha x^2-(\alpha^2+2\alpha-1)x+2\alpha^2-\alpha=0$. Moreover, by direct calculation, we have $\alpha>2-\frac{1}{\alpha}$ when $\alpha\in[\frac{1}{2},1)$. Together with \eqref{5.3}, we obtain
\begin{equation*}
x_{\max}^{2}\geq\alpha( t +2)+\frac{(1-\alpha)^2}{\alpha}>2\alpha>\alpha>2-\frac{1}{\alpha}.
\end{equation*}

Therefore, when $x= x_{\max}^{2}$, we have $f_1(\alpha, t ,x)-f_2(\alpha, t ,x)<0$, that is, $f_1(\alpha, t ,x)<0$. Notice that $f_1(\alpha, t ,x)$ is a monic polynomial, then $f_1(\alpha, t ,x)\geq0$ for $x\geq x_{\max}^{1}=\lambda_\alpha(T_1)$. Consequently, $\lambda_\alpha(T_2)<\lambda_\alpha(T_1)$.

Thus, the proof is complete.
\end{proof}

\begin{lem}
Let $T_1=T_2( t +1, t -1, t +1, t +2)$, $T_2=T_2( t +1, t +1, t -1, t +2)$. It follows that $\lambda_\alpha(T_1)>\lambda_\alpha(T_2)$ for $\frac{1}{2}\leqslant \alpha <1$.
\end{lem}
\begin{proof}
Let $W_1=N_{T_1}(u_1)\setminus \{w_1\}$, $W_2=N_{T_1}(u_2)\setminus \{w_1, w_2\}$, $W_3=N_{T_1}(u_3)\setminus \{w_2, w_3\}$, $W_4=N_{T_1}(u_4)\setminus\{w_3\}$. It is routine to check that $\pi_1=\{W_1\}\cup\{u_1\}\cup\{w_1\}\cup\{W_2\}\cup\{u_2\}\cup\{w_2\}\cup\{W_3\}\cup\{u_3\}\cup\{w_3\}\cup\{W_4\}\cup\{u_4\}$
is an equitable partition of $V(T_1)$, whose quotient matrix is given as $A_\alpha(T_1)_{\pi_1}$. 

Let $W_1'=N_{T_2}(u_1)\setminus \{w_1\}$, $W_2'=N_{T_2}(u_2)\setminus \{w_1, w_2\}$, $W_3'=N_{T_2}(u_3)\setminus \{w_2, w_3\}$, $W_4'=N_{T_2}(u_4)\setminus\{w_3\}$. Clearly, $\pi_2=\{W_1'\}\cup\{u_1\}\cup\{w_1\}\cup\{W_2'\}\cup\{u_2\}\cup\{w_2\}\cup\{W_3'\}\cup\{u_3\}\cup\{w_3\}\cup\{W_4'\}\cup\{u_4\}$
is an equitable partition of $V(T_2)$, whose quotient matrix is given as $A_\alpha(T_2)_{\pi_2}$. 

Let $f_1(\alpha, t ,x)=\det(xI-A_\alpha(T_1)_{\pi_1})$, $f_2(\alpha, t ,x)=\det(xI-A_\alpha(T_2)_{\pi_2})$. Assume that $x_{\max}^{1}$ (resp. $x_{\max}^{2}$) is the largest root of $f_1(\alpha, t ,x)=0$ (resp. $f_2(\alpha, t ,x)=0$). By the same discussion as that of Lemma \ref{lem5.6}, we have $x_{\max}^{1}=\lambda_\alpha(T_1)$ and $x_{\max}^{2}=\lambda_\alpha(T_2)$.
Together with Lemma \ref{lem4.2}, we have
\begin{equation}\label{5.4}
  \alpha( t +3)+\frac{(1-\alpha)^2}{\alpha}\leq\lambda_\alpha(T_2)=x_{\max}^{2}.
\end{equation}

Utilizing calculations by MATLAB R2023\cite{15}, we obtain
\begin{equation*}
  f_1(\alpha, t ,x)-f_2(\alpha, t ,x)=2(\alpha-1)^4(2\alpha-x)[\alpha x^2-(\alpha^2+2\alpha-1)x+2\alpha^2-\alpha]^2.
\end{equation*}

It is routine to check that $x_1:=\alpha$, $x_2:=2-\frac{1}{\alpha}$ are the roots of $\alpha x^2-(\alpha^2+2\alpha-1)x+2\alpha^2-\alpha=0$. By direct calculations, we have $\alpha>2-\frac{1}{\alpha}$ when $\alpha\in[\frac{1}{2},1)$. Together with \eqref{5.4}, we obtain
\begin{equation*}
x_{\max}^{2}\geq\alpha( t +3)+\frac{(1-\alpha)^2}{\alpha}>2\alpha>\alpha>2-\frac{1}{\alpha}.
\end{equation*}

Therefore, when $x= x_{\max}^{2}$, we have $f_1(\alpha, t ,x)-f_2(\alpha, t ,x)<0$, that is, $f_1(\alpha, t ,x)<0$. Notice that $f_1(\alpha, t ,x)$ is a monic polynomial, then $f_1(\alpha, t ,x)\geq0$ for $x\geq x_{\max}^{1}=\lambda_\alpha(T_1)$. Consequently, $\lambda_\alpha(T_2)<\lambda_\alpha(T_1)$.
\end{proof}

\begin{table}\caption{All the possible $A_\alpha$-minimizer graphs in $\mathcal{G}_{n,n-4}$ when $ t \geq3$.}\label{Tab 2}
  \centering
  \renewcommand{\arraystretch}{1.5}
  \begin{tabular}{p{1cm}<{\centering}|p{5cm}<{\centering}| p{5cm}<{\centering}}
\hline
$\ell'$&~$(\ell(w_1),\ell(w_2),\ell(w_3),\ell(w_4))$&$(\ell(u_1),\ell(u_2),\ell(u_3),\ell(u_4))$\\
\hline
0&\makecell*[l]{
{$( t ,~ t +1,~ t +1,~ t -2)$}\\
{$( t ,~ t ,~ t +1,~ t -1)$}\\
{$( t +1,~ t +1,~ t +1,~ t -3)$}}&
\makecell*[l]{
{$( t ,~ t ,~ t -1,~ t +1)$}\\
{$( t +1,~ t -1,~ t -1,~ t +1)$}}\\
\hline
1&\makecell*[l]{{$( t ,~ t +1,~ t +1,~ t -1)$}\\
{$( t +1,~ t +1,~ t +1,~ t -2)$}}&
\makecell*[l]{{$( t ,~ t ,~ t ,~ t +1)$}\\
{$( t +1,~ t -1,~ t ,~ t +1)$}}\\
\hline
2&\makecell[l]{{$( t +1,~ t +1,~ t +1,~ t -1)$}}&
\makecell[l]{{$( t +1,~ t ,~ t ,~ t +1)$}}\\
\hline
3&\makecell*[l]{
{$( t +1,~ t +2,~ t +2,~ t -2)$}\\
{$( t +1,~ t +1,~ t +2,~ t -1)$}\\
{$( t +1,~ t +1,~ t +1,~ t )$}\\
{$( t +2,~ t +2,~ t +2,~ t -3)$}}&
\makecell*[l]{
{$( t +1,~ t +1,~ t -1,~ t +2)$}\\
{$( t +1,~ t ,~ t ,~ t +2)$}\\
{$( t +2,~ t -1,~ t ,~ t +2)$}}\\
\hline
\end{tabular}
\end{table}

\begin{lem}\label{lem5.8}
Let $T_1=T_2( t +1, t , t +1, t +1)$, $T_2=T_2( t +1, t , t , t +2)$. Then one has $\lambda_\alpha(T_1)>\lambda_\alpha(T_2)$ for $\alpha\in[\frac{1}{2},1)$.
\end{lem}
\begin{proof}
Let $W_1=N_{T_1}(u_1)\setminus \{w_1\}$, $W_2=N_{T_1}(u_2)\setminus \{w_1, w_2\}$, $W_3=N_{T_1}(u_3)\setminus \{w_2, w_3\}$, $W_4=N_{T_1}(u_4)\setminus\{w_3\}$. It is routine to check that $\pi_1=\{W_1\}\cup\{u_1\}\cup\{w_1\}\cup\{W_2\}\cup\{u_2\}\cup\{w_2\}\cup\{W_3\}\cup\{u_3\}\cup\{w_3\}\cup\{W_4\}\cup\{u_4\}$
is an equitable partition of $V(T_1)$, whose quotient matrix is written as $A_\alpha(T_1)_{\pi_1}$. 

Let $W_1'=N_{T_2}(u_1)\setminus \{w_1\}$, $W_2'=N_{T_2}(u_2)\setminus \{w_1, w_2\}$, $W_3'=N_{T_2}(u_3)\setminus \{w_2, w_3\}$, $W_4'=N_{T_2}(u_4)\setminus\{w_3\}$. Clearly, $\pi_2=\{W_1'\}\cup\{u_1\}\cup\{w_1\}\cup\{W_2'\}\cup\{u_2\}\cup\{w_2\}\cup\{W_3'\}\cup\{u_3\}\cup\{w_3\}\cup\{W_4'\}\cup\{u_4\}$
is an equitable partition for $V(T_2)$, whose quotient matrix is given as $A_\alpha(T_2)_{\pi_2}$. 

Let $f_1(\alpha, t ,x)=\det(xI-A_\alpha(T_1)_{\pi_1})$, $f_2(\alpha, t ,x)=\det(xI-A_\alpha(T_2)_{\pi_2})$. Assume that $x_{\max}^{1}$ (resp. $x_{\max}^{2}$) is the largest zero of $f_1(\alpha, t ,x)$ (resp. $f_2(\alpha, t ,x)$). By the same discussion as that of Lemma \ref{lem5.6}, we have $x_{\max}^{1}=\lambda_\alpha(T_1)$ and $x_{\max}^{2}=\lambda_\alpha(T_2)$.
Together with Lemma \ref{lem4.2}, we have
\begin{equation}\label{5.5}
  \alpha( t +3)+\frac{(1-\alpha)^2}{\alpha}\leq\lambda_\alpha(T_2)=x_{\max}^{2}.
\end{equation}

Utilizing calculations by MATLAB R2023\cite{15}, we obtain 
$$f_1(\alpha, t ,x)-f_2(\alpha, t ,x)=(\alpha-1)^2 (2\alpha-x)(\alpha x-2\alpha+1)^2 g(x),$$

where
\begin{equation*}
\begin{split}
g(x) = &x^4 -(7\alpha+2\alpha t )x^3 +(\alpha^2{ t }\,^{2} +9\alpha^2  t +16\alpha^2 +2\alpha t +4\alpha- t -2)x^2\\
&+(-2\alpha^3 t \,^2 -10\alpha^3  t -14\alpha^3 -2\alpha^2 t \,^2 -12\alpha^2  t -12\alpha^2 +\alpha t \,^2 +6\alpha t +6\alpha)x\\
&+2\alpha^4  t +4\alpha^4 +4\alpha^3  t \,^2 +12\alpha^3  t +8\alpha^3 - 2\alpha^2 t \,^2-6\alpha^2 t -4\alpha^2,\\
g'(x) =&4x^3 -(21\alpha+6\alpha t )x^2 +(2\alpha^2  t \,^2 +18\alpha^2  t +32\alpha^2 +4\alpha t +8\alpha-2 t -4)x\\
&-2\alpha^3  t \,^2 -10\alpha^3  t -14\alpha^3 -2\alpha^2  t \,^2 -12\alpha^2  t -12\alpha^2 +\alpha t \,^2 +6\alpha t +6\alpha,\\
g''(x) =&12x^2 -(42\alpha+12\alpha t )x+2\alpha^2  t \,^2 +18\alpha^2  t +32\alpha^2 +4\alpha t +8\alpha-2 t -4.
\end{split}
\end{equation*}

Notice that $\alpha( t +3)+\frac{(1-\alpha)^2}{\alpha}>\frac{1}{2}\alpha t +\frac{7}{4}\alpha$. Thus, when $x\in[\alpha( t +3)+\frac{(1-\alpha)^2}{\alpha},+\infty)$, we have
$$g''(x)\geq g''(\alpha( t +3)+\frac{(1-\alpha)^2}{\alpha})=2\alpha^2  t \,^2 +(24\alpha^2-20\alpha+10) t +56\alpha^2-100\alpha+98-\frac{48}{\alpha}+\frac{12}{\alpha^2}>0.$$
Therefore, \(g'(x)\) behaves as an increasing function in $x$ for $x\in[\alpha( t +3)+\frac{(1-\alpha)^2}{\alpha},+\infty)$. Thus,
\begin{equation*}
\begin{split}
g'(x)\geq g'(\alpha( t +3)+\frac{(1-\alpha)^2}{\alpha})=&(3\alpha^3-2\alpha^2+\alpha) t \,^2 +(22\alpha^3-36\alpha^2+34\alpha-16+\frac{4}{\alpha}) t \\
&+34\alpha^3-92\alpha^2+138\alpha-124+\frac{71}{\alpha}-\frac{24}{\alpha^2}+\frac{4}{\alpha^3}>0,
\end{split}
\end{equation*}
which implies $g(x)$ behaves as an increasing function in $x$ for $x\in[\alpha( t +3)+\frac{(1-\alpha)^2}{\alpha},+\infty)$. Thus,
\begin{equation*}
\begin{split}
g(x)\geq g(\alpha( t +3)+\frac{(1-\alpha)^2}{\alpha})=&(2\alpha^4-2\alpha^3+\alpha^2) t \,^2 +\frac{(2\alpha^2-2\alpha+1)(6\alpha^4-8\alpha^3+8\alpha^2-4\alpha+1)}{\alpha^2}t \\
&+\frac{(2\alpha^2-2\alpha+1)(6\alpha^6-16\alpha^5+28\alpha^4-28\alpha^3+17\alpha^2-6\alpha+1)}{\alpha^4}>0.
\end{split}
\end{equation*}

It is straightforward to check that $x_1:=2\alpha$, $x_2=2-\frac{1}{\alpha}$ are two zeros of $f_1(\alpha, t ,x)-f_2(\alpha, t ,x)$. Moreover, by direct calculation, we have $2\alpha>2-\frac{1}{\alpha}$ when $\alpha\in[\frac{1}{2},1)$. Together with \eqref{5.5}, we obtain
\begin{equation*}
x_{\max}^{2}\geq\alpha( t +3)+\frac{(1-\alpha)^2}{\alpha}>2\alpha>2-\frac{1}{\alpha}.
\end{equation*}

Therefore, when $x= x_{\max}^{2}$, we have $f_1(\alpha, t ,x)-f_2(\alpha, t ,x)<0$, that is, $f_1(\alpha, t ,x)<0$. Notice that $f_1(\alpha, t ,x)$ is a monic polynomial, then $f_1(\alpha, t ,x)\geq0$ for $x\geq x_{\max}^{1}=\lambda_\alpha(T_1)$. Consequently, $\lambda_\alpha(T_2)<\lambda_\alpha(T_1)$.
\end{proof}

\begin{table}[h!]
  \caption{All the possible $A_\alpha$-minimizer graphs in $\mathcal{G}_{n,n-4}$ when $t=1,2$.}\label{Tab 3}
  \centering
  \renewcommand{\arraystretch}{1.5}
  \begin{tabular}{p{1cm}<{\centering}|p{1cm}<{\centering}|p{1cm}<{\centering}|p{4.5cm}<{\centering}| p{4.5cm}<{\centering}}
\hline
$t$&~$\ell'$&~$n$&~$(\ell(w_1),\ell(w_2),\ell(w_3),\ell(w_4))$&$(\ell(u_1),\ell(u_2),\ell(u_3),\ell(u_4))$\\
\hline
\multirow{3}{*}[-14pt]{1} &1&12&\makecell*[c]{
{$-$}}&
\makecell*[c]{
{$(2,~0,~1,~2)$}}\\
\cline{2-5}
&2&13&\makecell*[c]{{$(2,~2,~2,~0)$}}&
\makecell*[c]{$(2,~1,~1,~2)$}\\
\cline{2-5}
&3&14&\makecell*[c]{
{$(2,~2,~3,~0)$}\\
{$(2,~2,~2,~1)$}}&
\makecell*[c]{
{$(2,~2,~0,~3)$}\\
{$(2,~1,~1,~3)$}\\
{$(3,~0,~1,~3)$}}\\
\hline
\multirow{4}{*}[-20pt]{2} &0&15&\makecell*[c]{
{$(2,~3,~3,~0)$}\\
{$(2,~2,~3,~1)$}}&
\makecell*[c]{
{$(2,~2,~1,~3)$}\\
{$(3,~1,~1,~3)$}}\\
\cline{2-5}
&1&16&\makecell*[c]{
{$(2,~3,~3,~1)$}\\
{$(3,~3,~3,~0)$}}&
\makecell*[c]{
{$(2,~2,~2,~3)$}\\
{$(3,~1,~2,~3)$}}\\
\cline{2-5}
&2&17&\makecell*[c]{
{$(3,~3,~3,~1)$}}&
\makecell*[c]{
{$(3,~2,~2,~3)$}}\\
\cline{2-5}
&3&18&\makecell*[c]{
{$(3,~4,~4,~0)$}\\
{$(3,~3,~4,~1)$}\\
{$(3,~3,~3,~2)$}}&
\makecell*[c]{
{$(3,~3,~1,~4)$}\\
{$(3,~2,~2,~4)$}\\
{$(4,~1,~2,~4)$}}\\
\hline
\end{tabular}
\end{table}

\begin{remark}\label{remark2}
{\rm 
By above analysis and Lemmas \ref{lem5.4} and \ref{lem5.6}-\ref{lem5.8}, we can easily exclude some graphs from Table \ref{Tab 1} and obtain Table \ref{Tab 2}, which shows all the $18$ possible $A_\alpha$-minimizer graphs in $\mathcal{G}_{n,n-4}$ for $ t \geq3$. 

When $ t =1,2$, we obtain $12\leq n\leq 18$. As similar as that of case $ t \geq 3$, we can also exclude some graphs and list all the possible $A_\alpha$-minimizer graphs in Table \ref{Tab 3}. Moreover, by Theorem \ref{thm3.02}, one sees $W_{11}$ is the unique graph in $\mathcal{G}_{n,n-4}$ having minimal $A_\alpha$ spectral radius when $n=11$ and $\alpha\in[0,1)$.

Liu and Wang \cite{13} characterized the $A_\alpha$-minimizer graphs in $\mathcal{G}_{n,n-4}$ for $n\in\{5,7,8,9,10\}$ and $\alpha\in[0,1)$. When $n=6,$ they also identified the $A_\alpha$-minimizer graphs in $\mathcal{G}_{6,2}$ for $\alpha\in[0,\frac{7}{9}]$.  
}
\end{remark}
In the remainder of this section, we are to determine the $A_\alpha$-minimizer graphs in $\mathcal{G}_{6,2}$ for $\alpha\in[0,1)$. The next structural lemma is crucial to our proof for Theorem \ref{thm1.5}.

\begin{figure}[!ht]
\centering
\begin{tikzpicture}[scale = 1.2]
  \tikzstyle{vertex}=[circle,fill=black,minimum size=0.38em,inner sep=0pt]
  \node[vertex] (G_0) at (0,-0.2)[label=left:$v_1$]{};
  \node[vertex] (G_1) at (0.8,0.6)[label=left:$v_2$]{};
  \node[vertex] (G_2) at (0.8,-1)[label=left:$v_4$]{};
  \node[vertex] (G_3) at (0.8,-0.2)[label=above:$v_3$]{};
  \node[vertex] (G_4) at (1.6,0.2)[label=right:$v_i$]{};
  \node[vertex] (G_5) at (1.6,-0.8)[label=right:$v_j$]{};
  \draw[thick] (G_0) -- (G_1) -- (G_4) -- (G_3) -- (G_0);
  \draw[thick] (G_0) -- (G_2) -- (G_5) -- (G_4);
   \draw (0.8,-1.8)node{$G_1$};
  \end{tikzpicture}
   \hspace{2em}
   \begin{tikzpicture}[scale = 1.2]
  \tikzstyle{vertex}=[circle,fill=black,minimum size=0.38em,inner sep=0pt]
 \node[vertex] (G_0) at (0,-0.1)[label=above:$v_2$]{};
 \node[vertex] (G_1) at (-0.8,-0.4)[label=left:$v_4$]{};
 \node[vertex] (G_2) at (0.8,-0.4)[label=right:$v_j$]{};
 \node[vertex] (G_3) at (-0.5,0.4)[label=left:$v_i$]{};
 \node[vertex] (G_4) at (0.5,0.4)[label=right:$v_1$]{};
 \node[vertex] (G_5) at (0,-1)[label=below:$v_3$]{};
 \draw[thick] (G_3) -- (G_4) -- (G_2) -- (G_5) -- (G_1) -- (G_3);
 \draw[thick] (G_2) -- (G_0) -- (G_1);
 \draw[thick] (G_0) -- (G_5);
  \draw (0,-1.8)node{$G_2$};
  \end{tikzpicture}
  \hspace{2em}
\begin{tikzpicture}[thick]
   \tikzstyle{vertex}=[circle,fill=black,minimum size=0.38em,inner sep=0pt]
    \node[draw, circle, minimum size=1.5cm] (K_s) at (-1.5,-0.2) {$K_s$};
    \node[vertex, label={below right:$u$}] (u) at (-0.75,-0.2) {};
    \node[draw, circle, minimum size=1.5cm] (K_t) at (1.5,-0.2) {$K_t$};
    \node[vertex, label={below left:$v$}] (v) at (0.75,-0.2) {};
    \draw (u) -- (v);
     \draw (0,-1.8)node{$F_{s,t}$};
\end{tikzpicture}
    \caption{Graphs $G_1, G_2, F_{s,t}$ used in the proof of Lemma \ref{lem5.9} and Theorem \ref{thm1.5}.}\label{Fig.5}
\end{figure}
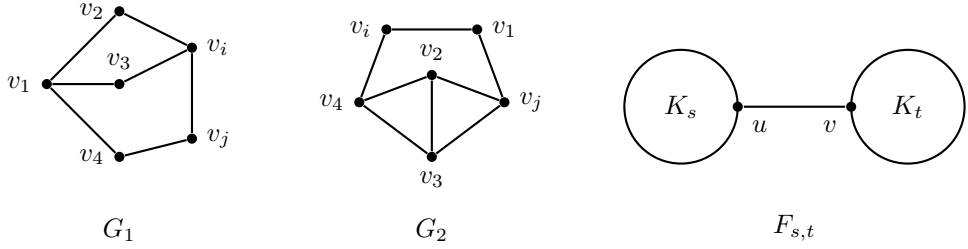

\begin{lem}\label{lem5.9}
If $G$ is in $\mathcal{G}_{6,2}$ having minimal \(A_\alpha\) spectral radius with $0\leqslant \alpha<1$, then $G^c$ is bipartite. 
\end{lem}
\begin{proof}
Note that $\gamma(G)=2$. Consequently, $G^c$ contains no $K_3$.  Select \(v_1 \in V(G)\) such that \(d_{G^c}(v_1) = \Delta(G^c)\). Given that \(G^c\) is \(K_3\)-free, it’s easy to verify that the induced subgraph \(G^c[N_{G^c}(v_1)]\) is an empty graph. To finish the proof, we just need to show that \(G^c[V(G) \setminus N_{G^c}(v_1)]\) is empty.

In light of Lemma \ref{lem5.9.1}, one sees  
$|E(G^c)|\leq\lfloor\frac{6^2}{4}\rfloor-1=8$, namely $|E(G)|\geq
\frac{6(6-1)}{2}-8$=7. We proceed by determine the range of $|E(G)|$.

In fact, if $|E(G)|\geq 9$, then by virtue of Lemma \ref{lem4.2.1}, one has 
\begin{equation}\label{5.6}
 \lambda_\alpha(G)\geq \frac{2|E(G)|}{|V(G)|}\geq3.
\end{equation}
It is clear that both equalities hold precisely when $G$ is regular and \(|E(G)| = 9\). It is straightforward to check that $F_{3,3}\in \mathcal{G}_{6,2}$ (see Fig. \ref{Fig.5}). Moreover, according to Lemma \ref{lem4.2.1}, one has $\lambda_\alpha(F_{3,3})\leq 3(1-\alpha)+3\alpha=3$. Therefore, we have 
\begin{equation}\label{5.7}
 \lambda_\alpha(G)\leq\lambda_\alpha(F_{3,3})\leq 3.
\end{equation}
Combining with \eqref{5.6} and \eqref{5.7} gives us $\lambda_\alpha(G)=3$ and $G$ is 3-regular with $|E(G)|=9$. Thus, $G$ is isomorphic to $K_{3,3}$ or the three-prism graph. Notice that $\gamma(K_{3,3})=3$, which leads to a contradiction since $\gamma(G)=2$. 
Moreover, it is straightforward to verify that \(F_{3,3}\) is properly contained in the three-prism graph, a contradiction to the minimality of $\lambda_\alpha(G)$.

Therefore, $|E(G)|=7~\mbox{or}~8$. Assume that $G^c[V(G)\setminus N_{G^c}(v_1)]$ is non-empty. It follows that there exist two vertices $v_i,v_j$ in $V(G)\setminus N_{G^c}(v_1)$ such that they are adjacent in $G^c$. Bearing in mind that $G^c$ contains no $K_3$, the number of edges formed by one vertex in $\{v_i,v_j\}$ and the other one in $N_{G^c}(v_1)$ is not more than $\Delta(G^c)$. Define $W:=V(G)\setminus(N_{G^c}(v_1)\cup\{v_i,v_j\})$. Then we have
\begin{align*}
   |E(G^c)|\leq\sum_{v\in W}d_{G^c}(v)+1+\Delta(G^c)
   \leq\Delta(G^c)(6-\Delta(G^c)-2)+1+\Delta(G^c)
   =-{\Delta(G^c)}^2+5\Delta(G^c)+1\leq7.
\end{align*}
Thus, $|E(G^c)|=-{\Delta(G^c)}^2+5\Delta(G^c)+1=7$, which implies that $\Delta(G^c)=2~\mbox{or}~3$.

If $\Delta(G^c)=2$, then $\delta(G)=3,$ and so $|E(G)|\geq \frac{|V(G)||\delta(G)|}{2}=9$, a contradiction.

If $\Delta(G^c)=3$, then $G^c\cong G_1$ and $G\cong G_{2}$ (see Fig. \ref{Fig.5}). 
Notice that $\pi_1=\{v_1,v_i\}\cup\{v_4,v_j\}\cup\{v_2,v_3\}$ is an equitable partition of $V(G_2)$, whose quotient matrix is $A_\alpha(G_2)_{\pi_1}$, $\pi_2=\{u,v\}\cup({V(F_{3,3})\setminus\{u,v\}})$ is an equitable partition of $V(F_{3,3})$, whose quotient matrix is $A_\alpha(F_{3,3})_{\pi_2}$.
\begin{equation*}
\renewcommand{\arraystretch}{1.3}
	A_\alpha(G_2)_{\pi_1}=\begin{pmatrix}
	1+\alpha&1-\alpha&0\\
	1-\alpha&3\alpha&2-2\alpha\\
	0&2-2\alpha&1+2\alpha
	\end{pmatrix},
    \hspace{2em}
    A_\alpha(F_{3,3})_{\pi_2}=\begin{pmatrix}
	1+2\alpha&2-2\alpha\\
	1-\alpha&1+\alpha
	\end{pmatrix}.
	\end{equation*}

Let $f_1(\alpha,x)=\det(xI-A_\alpha(G_2)_{\pi_1})$, $f_2(\alpha,x)=\det(xI-A_\alpha(F_{3,3})_{\pi_2})$. Assume that $x_{\max}^{1}$ (resp. $x_{\max}^{2}$) is the largest zero of  $f_1(\alpha,x)$ (resp. $f_2(\alpha,x)$). According to Lemma \ref{lem5.3}, we have $x_{\max}^{1}=\lambda_\alpha(G_2)$ and $x_{\max}^{2}=\lambda_\alpha(F_{3,3})$.

By direct calculations, we obtain
\begin{align*}
f_1(\alpha,x)&=x^3 -(6\alpha+2)x^2 +(6\alpha^2 +19\alpha-4)x-16\alpha^2 -7\alpha+5,\\
f_2(\alpha,x)&=x^2-(3\alpha+2)x+7\alpha-1,\\
x_{\max}^{2}&=\lambda_\alpha(F_{3,3})=\frac{3\alpha}{2}+\frac{\sqrt{9\alpha^2 -16\alpha+8}}{2}+1.\\
f_1(\alpha,x_{\max}^{2})&=-\frac{1}{2} (\alpha-1)^2 (3 \sqrt{9 \alpha^2-16 \alpha+8}+9 \alpha-4)<0.
\end{align*}
 Notice that $f_1(\alpha,x)$ is a monic polynomial. Hence, when $x\geq x_{\max}^{1}=\lambda_\alpha(G_2)$, one has $f_1(\alpha,x)\geq0$. Consequently, $\lambda_\alpha(F_{3,3})<\lambda_\alpha(G_2)$, a contradiction.
\end{proof}

\begin{thm}\label{thm1.5}
$F_{3,3}$ is the unique graph among $\mathcal{G}_{6,2}$ having minimal \(A_{\alpha}\) spectral radius, where $0\leqslant \alpha<1$. 
\end{thm}
\begin{proof}
By Lemmas \ref{lem2.1} and \ref{lem5.9}, we have $G\cong F_{3,3}, F_{2,4}~\mbox{or}~F_{1,5}$. Notice that $F_{2,4}$ is properly contained in $F_{1,5}$, again by virtue of Lemma \ref{lem2.1}, one has $\lambda_\alpha(F_{2,4})<\lambda_\alpha (F_{1,5})$. Subsequently, we proceed by proving $\lambda_\alpha(G_2)<\lambda_\alpha(F_{2,4})$.

Consider the graph $G_2$ (see Fig. \ref{Fig.5}). By virtue of Lemma \ref{lem5.9.2}, we obtain
\begin{align*}
   \lambda_\alpha(G_2)x_{v_i}=&2\alpha x_{v_i}+(1-\alpha)(x_{v_1}+x_{v_4})
   =2\alpha x_{v_i}+(1-\alpha)(x_{v_i}+x_{v_j}).
\end{align*}
Thus, we have
$(\lambda_\alpha(G_2)-\alpha-1)x_{v_i}=(1-\alpha)x_{v_j}.$
Acccording to Lemma \ref{lem4.2.1}, we have $\frac{8}{3}\leq\lambda_\alpha(G_2)\leq 3$. Then we obtain $x_{v_i}<x_{v_j}$. In light of Lemma \ref{lem2.2}, one has 
$$\lambda_\alpha(G_2)<\lambda_\alpha(G_2-v_iv_4+v_jv_4)=\lambda_\alpha(F_{2,4}).$$
Then by the proof of Lemma \ref{lem5.9}, we have $\lambda_\alpha (F_{3,3})<\lambda_\alpha(G_2)<\lambda_\alpha (F_{2,4})$. Therefore, the $A_\alpha$-minimizer graph in $\mathcal{G}_{6,2}$ is $F_{3,3}$ for $\alpha\in[0,1)$.
\end{proof}

By Remark~\ref{remark2} and Theorem~\ref{thm1.5}, all the possible $A_\alpha$-minimizer graphs in $\mathcal{G}_{n,n-4}$ for $n\geq5$ and $\alpha\in[\frac{1}{2},1)$ have been characterized.

\section*{\normalsize Disclosure statement}
The authors did not report any potential conflict of interest.

\section*{\normalsize Acknowledgments}
Shuchao Li receives financial support from the National Natural Science Foundation of China (Grant Nos. 12571365, 12171190), the Special Fund for Basic Scientific Research of Central Colleges (Grant Nos. CCNU25-\linebreak JC006, CCNU25HD044, CCNU25JCPT031) and the Open Research Fund of Key Laboratory of Nonlinear Analysis \& Applications (CCNU), Ministry of Education of China (Grant No. NAA2025ORG010).

\section*{\normalsize Data availability}
No data is available during the current study.


\begin{thebibliography}{99}
\small \setlength{\itemsep}{-.2mm}
\bibitem{11} J.A. Bondy, U.S.R. Murty, Graph theory with applications, London: Macmillan, 1976.
\bibitem{17}D.M. Cvetkovi\'c, M. Doob, H. Sachs, Spectra of Graphs, third ed., Johan Ambrosius Barth Verlag, Heidelberg, 1995.
\bibitem{16} P. Csikv\'ari, Integral trees of arbitrarily large diameters, J. Algebraic Comb. 32 (2010) 371–377.
\bibitem{10} H.Y. Guo, B. Zhou, On the $\alpha$-spectral radius of graphs, Appl. Anal. Discrete Math. 14(2) 2020 431-458.
\bibitem{8} Y.R. Hu, Q.X. Huang, Z.Z. Lou, Graphs with the minimum spectral radius for given independence number, Discrete Math. 348 (2025) 114265.
\bibitem{12} Y.R. Hu, Z.Z. Lou, W.J. Ning, The $Q$-minimizer graph with given independence number, Linear Algebra Appl. 685 (2024) 1-23.
\bibitem{5} C.Y. Ji, M. Lu, On the spectral radius of trees with given independence number, Linear Algebra Appl. 488 (2016) 102-108.
\bibitem{LL2025}X.Y. Lei, S.C. Li, Spectral extremal results on the  $A_\alpha$-spectral radius of graphs without $K_{a,b}$-minor, 
Appl. Math. Comput. 492 (2025), Paper No. 129232, 21 pp. 
\bibitem{6} D. Li, Y.Y. Chen, J.X. Meng, The $A_\alpha$-spectral radius of trees and unicyclic graphs with given degree sequence, Appl. Math. Comput. 363 (2019) 124622.
\bibitem{13} X.C. Liu, L.G. Wang, The $A_\alpha$ spectral radius of graphs with given independence number $n-4$, Discrete Optim. 59 (2026) 100930.
\bibitem{7} Z.Z. Lou, J.M. Guo, The spectral radius of graphs with given independence number, Discrete Math. 345 (2022) 112778.
\bibitem{1} V. Nikiforov, Merging the $A$- and $Q$-spectral theories, Appl. Anal. Discrete Math. 11 (2017) 81-107.
\bibitem{3} V. Nikiforov, G. Pastn, O. Rojo, R.L. Soto, On the $A_\alpha$-spectra of trees, Linear Algebra Appl. 520 (2017) 286-305.
\bibitem{4} W.T. Sun, L.X. Yan, S.C. Li, X.C. Li, Sharp Bounds on the $A_\alpha$-index of Graphs In Terms of the Independence Number, Acta Math. Appl. Sin. 39 (2023) 656-674.
\bibitem{15} The MathWorks Inc., MATLAB, version R2023a, 2023, The MathWorks Inc. Natick.
\bibitem{9} M.M. Xu, Y. Hong, J.L. Shu, M.Q. Zhai, The minimum spectral radius of graphs with given independence number, Linear Algebra Appl. 431 (2009) 937-945.
\bibitem{2} J. Xue, H.Q. Lin, S.T. Liu, J.L. Shu, On the $A_\alpha$-spectral radius of a graph, Linear Algebra Appl. 550 (2018) 105-120.
\bibitem{14} L.H. You, M. Yang, W. So, W.G. Xi, On the spectrum of an equitable quotient matrix and its application, Linear Algebra Appl. 577 (2019) 21-40.
\bibitem{YGL2024}Y.T. Yu, X.Y. Geng, S.C. Li, An $A_\alpha$-spectral version of the Bhattacharya-Friedland-Peled conjecture, 
Electron. J. Combin. 31 (2024), no. 4, Paper No. 4.78, 27 pp.
\bibitem{ZL2025}Z.H. Zhou, S.C. Li, On the $A_\alpha$-index of graphs with given order and dissociation number, Discrete Appl. Math. 360 (2025), 167-180.
\end{thebibliography}
\end{document}